\documentclass{amsart}
\usepackage{amsmath,amssymb}

\def\normo#1{\left\|#1\right\|}

\def\brk#1{\left(#1\right)}

\def\rev#1{\frac{1}{#1}}
\def\half#1{\frac{#1}{2}}
\def\norm#1{\|#1\|}
\def\jb#1{\langle#1\rangle}
\def\wt#1{\widetilde{#1}}
\def\wh#1{\widehat{#1}}

\newcommand{\T}{{\mathbb T}}

\newcommand{\R}{{\mathbb R}}
\newcommand{\C}{{\mathbb C}}
\newcommand{\Z}{{\mathbb Z}}
\newcommand{\ft}{{\mathcal{F}}}

\newcommand{\les}{{\lesssim}}
\newcommand{\ges}{{\gtrsim}}

\newtheorem{theorem}{Theorem}[section]
\newtheorem{proposition}[theorem]{Proposition}
\newtheorem{lemma}[theorem]{Lemma}

\theoremstyle{definition}

\theoremstyle{remark}

\numberwithin{equation}{section}

\begin{document}

\title[Sharp Well-posednss for the Benjamin Equation]{Sharp Well-posedness for the
Benjamin Equation}


\thanks{The first/third-named author is supported in part by NNSF of China
 (No.10771130, No.10931001)/NESEC of Canada
(No. RGPIN/261100-2003).}

\author{Wengu Chen}
\address{Institute of Applied Physics
and Computational Mathematics, P.O.Box 8009, Beijing 100088, China}
\email{chenwg@iapcm.ac.cn}
\author{Zihua Guo}
\address{LMAM, School of Mathematical Sciences, Peking University,
Beijing 100871, China} \email{guozihua@gmail.com}
\author{Jie Xiao}
\address{Department of Mathematics and Statistics, Memorial University of Newfoundland, St. John's, NL A1C 5S7, Canada}
\email{jxiao@mun.ca}


\date{}



\begin{abstract}
     Having the ill-posedness in the range $s<-3/4$ of the Cauchy problem for the Benjamin equation with an initial $H^{s}({\mathbb R})$ data, we prove that the already-established local well-posedness in the range $s>-3/4$ of this initial value problem is extendable to $s=-3/4$ but also that such a well-posed property is globally valid for $s\in [-3/4,\infty)$.
\end{abstract}
\maketitle

\section{Introduction}

Continuing from \cite{CX}, we fully investigate the local and global well-posedness of the initial value problem for the Benjamin equation
 \begin{equation}\label{Ben1}
\left\{
\begin{array}{l}
\partial_t u - \gamma\partial_x u+\alpha {\mathcal
H}\partial^2_xu+\beta\partial^3_xu+\partial_x(u^2) =0, \quad
(x,t)\in {\mathbb R\times\mathbb
R},\\
u(x,0)=u_0(x), \quad x\in\mathbb R.
\end{array}
\right.
\end{equation}

Here ${\mathcal H}$ stands for the one-dimensional Hilbert transform:
$$
{\mathcal H}f(x)=\lim_{\epsilon\to 0}\frac1\pi\int_{|y|>\epsilon}
{f(x-y)}y^{-1}dy,\quad x\in \mathbb R.
$$
In Physics, the Benjamin equation, according to \cite{B} and \cite{Pe}, describes the vertical displacement, bounded above and below by rigid horizontal planes, of the interface between a thin layer of fluid atop and a much thicker layer of higher density fluid; see also \cite{ABR, A, B, CB, LS} for the study of
existence, stability and asymptotics of solitary wave solutions of (\ref{Ben1}). In addition, the case $\alpha\not=0$ and $\beta=0$ in (\ref{Ben1}) induces the Benjamin-Ono equation -- see
Kenig's survey \cite{Ke} but also Ionescu-Kenig \cite{IK} and
Burq-Planchon \cite{BP} for more information.

In Chen-Xiao's paper \cite{CX}, the following well/ill-posedness result was established through a sharp bilinear estimate for the so-called Bourgain space \cite{Bo} and Bejenaru-Tao's argument for \cite[Theorem 2]{BeTa} plus an example in Bourgain
\cite{Bo97} and Tzvetkov \cite{Tz}.

\begin{theorem}\label{thm2} For $\alpha,\beta,\gamma,\xi,s,b\in {\mathbb R}$ with
$\alpha\beta\not=0$, let $p(\xi)=\beta\xi^3-\alpha\xi|\xi|+\gamma\xi$. Then:

\item{\rm(i)} For $s> -3/4$ and $u_0\in H^s({\mathbb R})$ there exist
$b\in (1/2,1)$ and $T=T(\|u_0\|_{H^s})>0$ such that (\ref{Ben1}) has
a unique solution $u$ in $X_{s,b,p}\cap C([-T,T]; H^s({\mathbb R}))$.

\item{\rm(ii)} For $s<-3/4$ the solution map of the initial value problem of (\ref{Ben1}) is
not $C^3$ smooth at zero -- there is no $T>0$ such that the solution map of (\ref{Ben1}):
$$
u_0\in H^s({\mathbb R})\mapsto u\in C([-T,T]; H^s({\mathbb R}))
$$
is $C^3$ smooth at zero.
\end{theorem}

In the above and below, $H^s:=H^s(\mathbb R)$ and $X_{s,b,p}:=X_{s,b,p}(\mathbb R^2)$ represent the square Sobolev space with order $s$ and the Bourgain space -- the completion of all
$C_0^\infty(\mathbb R)$ and $C_0^\infty(\mathbb R^2)$ functions $f$ whose Fourier transforms $\widehat{f}$ satisfy
$$
\|f\|_{H^s}=\left(\int_\mathbb R
(1+|\xi|^2)^s|\widehat{f}(\xi)|^2\,d\xi\right)^\frac12<\infty
$$
and
$$
\|f\|_{X_{s,b,p}}=\left(\int_{\mathbb R}\int_{\mathbb
R}(1+|\xi|^2)^{s}\big(1+|\tau-p(\xi)|^2\big)^{b}\,
|\widehat{f}(\xi,\tau)|^2 \,d\xi\,d\tau\right)^{1/2}<\infty
$$
respectively. Moreover $H^\infty(\mathbb R):=\bigcap_{s=1}^\infty
H^s(\mathbb R)$ is equipped with the induced metric. Here it is appropriate to
point out that: Under $(\beta,\gamma)=(-1,0)$ and $\alpha=-\nu\in
(-1,0)$ Theorem \ref{thm2} (i) returns to Kozono-Ogawa-Tanisaka's
\cite[Theorem 2.1]{KOT}; Under $s\ge-1/8$ and $\gamma=0$ Theorem
\ref{thm2} (i) goes back to Guo-Huo's \cite[Theorem 1.1]{GH}; Under
$(s,\gamma)=(0,0)$ and $\alpha\beta>0$ Theorem \ref{thm2} (i) yields
the local well-posedness in Linares \cite{L}.

Due to the fact that the key bilinear estimate for
a Bourgain space (stated in \cite[Theorem 1.1]{CX}):
\begin{equation}\label{bil}
\|\partial_x(uv)\|_{X_{s,\,(b-1)+,p}}\le c
\|u\|_{X_{s,b,p}}\|v\|_{X_{s,b,p}}
\end{equation}
fails for any $s\le -3/4$ and $b\in {\mathbb R}$, the paper \cite{CX} conjectured that (\ref{Ben1}) is locally well-posed for the intermediate index $s=-3/4$. The first aim of this paper is to verify this conjecture by modifying $X_{s,b,p}$ in terms of an $l^1$ Besov-type norm -- such a modification is mainly motivated by: Tataru's \cite{Ta} on wave maps; Bejenaru-Tao's \cite{BeTa} on
Schr\"{o}dinger equation; Ionescu-Kenig's \cite{IK} on BO equation; Ionescu-Kenig-Tataru's \cite{IKT} on KP-I equation; Guo's \cite{Guo} on KdV equation.

\begin{theorem}\label{thm3} For $u_0\in H^{-3/4}(\mathbb R)$ there exists $T=T(\|u_0\|_{H^{-3/4}})$ such that the initial value problem (\ref{Ben1})
has a solution $u$ in $\bar{F}^{-3/4}\cap C([-T,T];
H^{-3/4}(\mathbb R))$, but also the solution map $u_0\mapsto u$ is
the unique extension of the classical solution map from
$H^\infty(\mathbb R)$ into $C([-T,T]; H^\infty(\mathbb R))$.
\end{theorem}

On the other hand, since Linares \cite{L} obtained the
global well-posedness for (\ref{Ben1}) at $s=0$ only via the $L^2$
conservation law, the second aim of this paper is to effectively
adapt both the I-method developed by
Colliander-Keel-Staffilani-Takaoka-Tao in \cite{I-method} and the
approach taken in Guo \cite{Guo} to show that the solutions in
Theorems \ref{thm2}-\ref{thm3} actually exist for $t$ in an
arbitrary time interval $[0,T]$, thereby establishing the sharp
global well-posedness of (\ref{Ben1}) below.

\begin{theorem}\label{thm4} The initial value problem (\ref{Ben1})
is globally well-posed for $u_0\in H^s(\mathbb R)$ with $s\in [-3/4,\infty)$.
\end{theorem}

Before verifying Theorems \ref{thm3}-\ref{thm4} in Sections 2-3-4-5, let us agree to several basic notations. As above, by the Fourier transform $\hat{f}$ (or ${\mathcal F}(f)$) of $f\in {\mathcal S}'(\mathbb R^2)$ we mean:
$$
\hat{f}(\xi,\,\tau)=\int_{{\mathbb R}}\int_{\mathbb R}e^{-i(t\tau+x\xi)}f(x,\,t)dxdt.
$$
For the integer set $\mathbb Z$, let ${\mathbb Z}_+=\mathbb Z\cap [0,\infty)$ and
\[
I_k=\left\{
\begin{array}{l}
\{\xi:\,|\xi|\in [2^{k-1},\,2^{k+1}]\}\quad\hbox{when}\quad 0<k\in\mathbb Z_+,\\
\{\xi:\,|\xi|\leq 2\}\quad\hbox{when}\quad k=0.
\end{array}
\right.
\]

Denote by $\eta_0:\,{\mathbb R}\rightarrow [0,\,1]$ a bump function
adapted to $[-8/5,\,8/5]$ and take value $1$ in $[-5/4,\,5/4]$. For $k\in {\mathbb Z}$ set
\[
\eta_k(\xi)\equiv
\left\{
\begin{array}{l}
\eta_0(\xi/2^k)-\eta_0(\xi/2^{k-1})\quad \hbox{for}\quad k\geq 1,\\
0\quad\hbox{for}\quad k\leq -1.
\end{array}
\right.
\]
For $k\in {\mathbb Z}$ let $\chi_k(\xi)=\eta_0(\xi/2^k)-\eta_0(\xi/2^{k-1})$. Following \cite{IKT}, given $k\in {\mathbb Z}_+$ define
\begin{equation*}
X_k=\{f\in L^2({\mathbb R}^2)\mbox{ with support in } I_k\times
{\mathbb R}\mbox{ such that } \|f\|_{X_k}<\infty\}
\end{equation*}
as the dyadic $X^{s,\,b}$ type space, where
\begin{equation*}
\|f\|_{X_k}=\sum_{j=0}^\infty 2^{j/2}\|\eta_j(\tau-p(\xi))\cdot
f\|_{L^2_{\xi,\,\tau}}.
\end{equation*}
The $l^1$-analogue $F^s$ of an $X^{s,\,b}$ space, as in
\cite{BeTa} and \cite{Guo}, is determined by
\begin{equation*}
\|u\|_{F^s}^2=\sum_{k\geq 0}2^{2sk}\|\eta_k(\xi){\mathcal
F}(u)\|_{X_k}^2.
\end{equation*}

Denote by $A\lesssim B$ the inequality that $A\leq CB$ holds for some large constant $C$ that
may change (line by line) and rely on various parameters;
similarly employ $A\ll B$ to represent $A\leq C^{-1}B$; use $A\sim
B$ to stand for $A\lesssim B\lesssim A$; and write $<\xi>=(1+|\xi|^2)^{1/2}$ when $\xi\in\mathbb R$. So, from the definition of $X_k$ we can see that for any $l\in{\mathbb
Z}_+$ and $f_k\in X_k$ (cf. \cite{IKT}),
\begin{equation*}
\sum_{j=0}^\infty
2^{j/2}\Big\|\eta_j(\tau-p(\xi))\int|f_k(\xi,\,\tau')|2^{-l}
(1+2^{-l}|\tau-\tau'|)^{-4}d\tau'\Big\|_{L^2}\lesssim\|f_k\|_{X_k}.
\end{equation*}
Consequently, for $l\in{\mathbb Z}_+,\,t_0\in{\mathbb R},\,f_k\in X_k$
and $\gamma\in {\mathcal S}({\mathbb R})$ we have
\begin{equation}\label{eq05}
\|{\mathcal F}[\gamma(2^l(t-t_0))\cdot {\mathcal
F}^{-1}f_k]\|_{X_k}\lesssim \|f_k\|_{X_k}.
\end{equation}
Under $k\in {\mathbb Z}$ let $P_k$ stand for the operator on
$L^2({\mathbb R})$ defined by
$$
\widehat{P_ku}(\xi)=\eta_k(\xi) \hat{u}(\xi).
$$
For our convenience, we take a slight abuse of notation that $P_k$ is also treated as an operator on
$L^2({\mathbb R}\times {\mathbb R})$ by the formula
$$
{\mathcal
F}(P_ku)(\xi,\,\tau)=\eta_k(\xi){\mathcal F}(u)(\xi,\,\tau).
$$
Naturally, for $l\in {\mathbb Z}$ we put
$$
P_{\leq l}=\sum_{k\leq l}P_k,\quad P_{\geq l}=\sum_{k\geq l}P_k.
$$
In order to avoid some logarithmic divergence, we need to use a
weaker norm for the low frequency as in \cite{Guo}
$$
\|u\|_{\bar{X}_0}=\|u\|_{L_x^2L_t^\infty}.
$$
When $-3/4\leq s\leq 0$, we define the normed spaces:
\begin{equation*}
\bar{F}^s=\Big\{u\in {\mathcal S}'({\mathbb
R}^2):\,\|u\|_{\bar{F}^s}^2=\sum_{k\geq
1}2^{2sk}\|\eta_k(\xi){\mathcal
F}(u)\|_{X_k}^2+\|P_{0}(u)\|_{\bar{X}_0}^2<\infty\Big\}.
\end{equation*}
And for each $T>0$, we define the time-localized space $\bar{F}^s(T)$ through
\begin{equation*}
\|u\|_{\bar{F}^s(T)}=\inf_{v\in
\bar{F}^s\ with\ v=u\ {on}\ [-T,\,T]}\big(\|P_0u\|_{L_x^2L_{t\in [-T,T]}^\infty}+\|P_{\geq
1}v\|_{\bar{F}^s}\big).
\end{equation*}

Other notations are introduced during the developments that come up in the subsequent sections

\section{Dyadic Estimates for Local Well-posedness}

In this section we present several dyadic estimates lemmas.

\begin{lemma}[estimates for free Benjamin equation]\label{lem1}
For $t\in{\mathbb R}$ let $W(t)$ denote the solution at time $t$ of
the free Benjamin evolution, i.e., the operator on $L^2({\mathbb
R})$ defined by the Fourier multiplier $e^{itp(\xi)}$. Suppose
$I\subset{\mathbb R}$ is a time interval with $|I|\lesssim 1$ and
$k\in {\mathbb Z}_+$ and $k\geq 10$. If $\phi\in {\mathcal
S}({\mathbb R})$, then:
\begin{equation}\label{free}
\left\{
\begin{array}{l}
\|W(t)P_k\phi\|_{L_t^qL_x^r}\lesssim 2^{-k/q}\|\phi\|_{L^2};\\
\|W(t)P_{k}\phi\|_{L_x^2L_{t\in I}^\infty}\lesssim
2^{3k/4}\|\phi\|_{L^2};\\
\|W(t)P_k\phi\|_{L_x^4L_t^\infty}\lesssim 2^{k/4}\|\phi\|_{L^2};\\
\|W(t)P_k\phi\|_{L_x^\infty L_t^2}\lesssim 2^{-k}\|\phi\|_{L^2},
\end{array}
\right.
\end{equation}
where $(q,\,r)$ satisfies $2\leq q,\,r\leq\infty$ and $2/q=1/2-1/r$.
\end{lemma}

\begin{proof} For the first inequality, see \cite{GPW}, for the second see
\cite{KPVJAMS91}. For the third we use the results in \cite{KPV91},
for the last we use the results in \cite{KPVJAMS91} by noting that
$|p'(\xi)|\sim 2^{2k}$ if $|\xi|\sim 2^{k}$.
\end{proof}

\begin{lemma}[$X_k$ embedding]\label{lem2} Suppose $I\subset{\mathbb R}$ is a time interval with $|I|\lesssim 1$ and
$k\in {\mathbb Z}_+$ and $k\geq 10$. Let $(q,\,r)$ be defined as in
Lemma \ref{lem1}. If $\phi\in {\mathcal S}({\mathbb R})$, then
\begin{equation}\label{xk}
\left\{
\begin{array}{l}
\|P_k(u)\|_{L_t^qL_x^r}\lesssim 2^{-k/q}\|{\mathcal F}[P_k(u)]\|_{X_k},\\
\|P_{k}(u)\|_{L_x^2L_{t\in I}^\infty}\lesssim 2^{3k/4}\|{\mathcal F}[P_k(u)]\|_{X_k},\\
\|P_k(u)\|_{L_x^4L_t^\infty}\lesssim 2^{k/4}\|{\mathcal F}[P_k(u)]\|_{X_k},\\
\|P_k(u)\|_{L_x^\infty L_t^2}\lesssim 2^{-k}\|{\mathcal
F}[P_k(u)]\|_{X_k}.
\end{array}
\right.
\end{equation}
Moreover, $u\in{\bar{F}}^s$ implies
$\|u\|_{L_t^\infty H^s}\lesssim \|u\|_{{\bar{F}}^s}.$
\end{lemma}

\begin{proof} It follows from (\ref{free}) in Lemma \ref{lem1} and a suitable adaption of \cite[Lemma 3.2]{Guo} for KdV equation.
\end{proof}

To see the next lemma, we need a few more definitions. For $k\in{\mathbb Z}$ and $j\in{\mathbb Z}_+$ we define
$$
D_{k,\,j}=\{(\xi,\,\tau):\,\xi\in
[2^{k-1},\,2^{k+1}]\, \mbox{and}\,\, \tau-p(\xi)\in I_j\}.
$$
For any
$k_1,\,k_2,\,k_3\in{\mathbb Z}$ and $j_1,\,j_2,\,j_3\in{\mathbb
Z}_+$, we consider
\begin{equation*}\label{block}
\|\chi\|_D:=\sup_{(u_{k_2,j_2},\,v_{k_3,j_3})\in
E}\|\chi_{D_{k_1,j_1}}(\xi,\,\tau)\cdot u_{k_2,j_2}*
v_{k_3,j_3}(\xi,\,\tau)\|_{L^2_{\xi,\tau}}
\end{equation*}
where the supremum is taken over
$$
E:=\Big\{(u,\,v):\,\|u\|_2,\,\|v\|_2\leq
1\,\,\mbox{and}\,\,\mbox{supp}(u)\subset
D_{k_2,j_2},\,\,\mbox{supp}(v)\subset D_{k_3,j_3}\Big\}.
$$
At the same time, we recall some of Tao's notations in \cite{Tao2001}. Any summations over capitalized variables such as $N_j,\,L_j,\,H$ are presumed to be dyadic, namely,
these variables range over numbers of the form $2^k$
for $k\in {\mathbb Z}$. The symbols $N_{max}, N_{med}, N_{min}$ stand for the
maximum, median, and minimum of three positive numbers $N_1,\,N_2,\,N_3$ respectively, and
hence $N_{max}\ge N_{med}\ge N_{min}$. Similarly, one has
$L_{max}\geq L_{med}\geq L_{min}$ when $L_1,\,L_2,\,L_3>0$. More
than that, we adopt the following summation convention: Any
summation of the form $L_{max}\sim\cdots$ is a sum over the three
dyadic variables $L_1,\,L_2,\,L_3\gtrsim 1$, for instance,
$$
\sum_{L_{max}\sim H}:=\sum_{L_1,\,L_2,\,L_3\gtrsim 1:\,L_{max}\sim
H}.
$$
Likewise, any summation of the form $N_{max}\sim\cdots$ sum over
the three dyadic variables $N_1,\,N_2,\,N_3>0$, in particular,
$$
\sum_{N_{max}\sim N_{med}\sim
N}:=\sum_{N_1,\,N_2,\,N_3>0:\,N_{max}\sim N_{med}\sim N}.
$$

So, it is easy to see that in order for $\|\chi\|_D$ to be nonzero, one must
require
\begin{equation}\label{double}
|k_{max}-k_{med}|\leq 3\quad\hbox{and}\quad 2^{j_{max}}\thicksim\max(2^{j_{med}},\,2^{2k_{max}+k_{min}}).
\end{equation}

\begin{lemma}[block estimates]\label{block estimate} Let $k_1,\,k_2,\,k_3\in{\mathbb Z}$
and $j_1,\,j_2,\,j_3\in{\mathbb Z}_+$ obey (\ref{double}). Suppose
$N_i=2^{k_i},\,L_i=2^{j_i}$ for $i=1,2,3$. Then:

\item{\rm(i)} $N_{max}\sim N_{min}\ \&\ L_{max}\sim N_{max}^2N_{min}$ implies
\begin{eqnarray}\label{estimate1}
\|\chi\|_D\lesssim L_{min}^{1/2}N_{max}^{-1/4}L_{med}^{1/4}.
\end{eqnarray}

\item{\rm(ii)} Anyone of the following three conditions
$$
\left\{
\begin{array}{l}
N_1\sim N_2\gg N_3\ \& \ N_{max}^2N_{min}\sim L_3\gtrsim L_2,\,L_1;\\
N_2\sim N_3\gg N_1\ \& \ N_{max}^2N_{min}\sim L_1\gtrsim L_2,\,L_3;\\
N_3\sim N_1\gg N_2\ \& \ N_{max}^2N_{min}\sim L_2\gtrsim L_3,\,L_1,
\end{array}
\right.
$$
implies
\begin{eqnarray}\label{estimate2}
\|\chi\|_D\lesssim
L_{min}^{1/2}N_{max}^{-1}\Big(\min\big\{N_{max}^2N_{min},\,
\frac{N_{max}}{N_{min}}L_{med}\big\}\Big)^{1/2}.
\end{eqnarray}

\item{\rm(iii)} In all other cases, one has
\begin{eqnarray}\label{estimate3}
\|\chi\|_D\lesssim
L_{min}^{1/2}N_{max}^{-1}\Big(\min\big\{N_{max}^2N_{min},\,L_{med}\big\}\Big)^{1/2}.
\end{eqnarray}
\end{lemma}
\begin{proof} It follows from \cite[Lemma 2.2]{CX}.
\end{proof}

Based on (\ref{estimate1})-(\ref{estimate2})-(\ref{estimate3}) of Lemma \ref{block estimate}, we obtain the forthcoming four dyadic bilinear bounds.

\begin{lemma}[high -- low interaction]\label{h-l}

\item{\rm (i)} If $k\geq 0,\,|k-k_2|\leq 5$, then for any $u,\,v\in
\bar{F}^s$,

$$
\Big\|(i+\tau-p(\xi))^{-1}\eta_k(\xi)i\xi\widehat{P_{0}u}*\widehat{P_{k_2}v}
\Big\|_{X_k}\lesssim\ \|P_{
0}u\|_{L_x^2L_t^\infty}\|\widehat{P_{k_2}v}\|_{X_{k_2}}.
$$

\item{\rm (ii)} If $k\geq 0,\,|k-k_2|\leq 5$ and $1\leq k_1\leq k-9$ then for any $u,\,v\in
\bar{F}^s$,
$$\Big\|(i+\tau-p(\xi))^{-1}\eta_k(\xi)i\xi\widehat{P_{k_1}u}*\widehat{P_{k_2}v}
\Big\|_{X_k}\lesssim\ k^32^{-k/2}2^{-k_1}\|\widehat{P_{
k_1}u}\|_{X_{k_1}}\|\widehat{P_{k_2}v}\|_{X_{k_2}}.
$$
\end{lemma}
\begin{proof} Without loss of generality, we may assume $k=k_2$.

(i) From the definition of $X_k$ it follows that
$$
\Big\|(i+\tau-p(\xi))^{-1}\eta_k(\xi)i\xi\widehat{P_{0}u}*\widehat{P_{k}v}
\Big\|_{X_k}\lesssim\ 2^{k}\sum_{j\geq 0}2^{-j/2}\|\widehat{P_{
0}u}*\widehat{P_{k}v}\|_{L^2_{\xi,\tau}}.
$$
By the Plancherel theorem
and (\ref{xk}) in Lemma \ref{lem2}, we get
$$2^{k}\|\widehat{P_{
0}u}*\widehat{P_{k}v}\|_{L^2_{\xi,\tau}}\lesssim
2^k\|P_0u\|_{L^2_xL_t^\infty}\|P_kv\|_{L_x^\infty L_t^2}\lesssim
\|P_0u\|_{L^2_xL^\infty_t}\|\widehat{P_kv}\|_{X_k}.
$$

(ii) Suppose
\begin{equation}\label{u-v}
u_{k_1,j_1}=\eta_{k_1}(\xi)\eta_{j_1}(\tau-p(\xi))\hat{u},\,v_{k,j_2}=
\eta_k(\xi)\eta_{j_2}(\tau-p(\xi))\hat{v}.
\end{equation}
Then
\begin{eqnarray}\label{v-u}
&&\Big\|(i+\tau-p(\xi))^{-1}\eta_k(\xi)i\xi\widehat{P_{k_1}u}*\widehat{P_{k}v}
\Big\|_{X_k}\nonumber\\
&&\lesssim 2^{k}\sum_{j_i\geq
0}2^{-j_3/2}\|\chi_{D_{k,j_3}}\cdot
u_{k_1,j_1}*v_{k,j_2}\|_{L^2_{\xi,\tau}}.
\end{eqnarray}
The estimate (\ref{double}) allows us to assume that $j_{max}\geq 2k+k_1-10$ in the summation on the right-hand side of (\ref{v-u}). Meanwhile we may also assume that $j_i\leq
10k\,(i=1,2,3)$ since otherwise an application of the trivial estimate
$$\|\chi_{D_{k,j_3}}\cdot
u_{k_1,j_1}*v_{k,j_2}\|_{L^2_{\xi,\tau}}\lesssim
2^{j_{min}/2}2^{k_{min}/2}\|u_{k_1,j_1}\|_{L^2_{\xi,\tau}}\|v_{k,j_2}\|_{L^2_{\xi,\tau}}
$$
gives the desired bound. Upon applying (\ref{estimate2}) we get
\begin{eqnarray*}
&& 2^k\sum_{j_i\geq 0}2^{-j_3/2}\|\chi_{D_{k,j_3}}\cdot
u_{k_1,j_1}*v_{k,j_2}\|_{L^2_{\xi,\tau}}\\
&&\lesssim 2^k\sum_{j_i\geq
0}2^{-j_3/2}2^{j_{min}/2}2^{-k/2}2^{-k_1/2}2^{j_{med}/2}\|u_{k_1,j_1}\|_{L^2_{\xi,\tau}}\|v_{k,j_2}\|_{L^2_{\xi,\tau}}\\
&&\lesssim 2^k\sum_{j_{max}\geq
2k+k_1-10}k^32^{-k/2}2^{-k_1/2}2^{-j_{max}/2}\|\widehat{P_{k_1}u}\|_{X_{k_1}}\|\widehat{P_kv}\|_{X_k}\\
&&\lesssim
k^32^{-k/2}2^{-k_1/2}\|\widehat{P_{k_1}u}\|_{X_{k_1}}\|\widehat{P_kv}\|_{X_k},
\end{eqnarray*}
thereby reaching the desired bound.
\end{proof}

When the low frequency is comparable to the high frequency, we have the following lemma.
\begin{lemma}[low $\sim$ high interaction]\label{h-l2} If $k\geq 10,\,|k-k_2|\leq 5$ and $k-9\leq k_1\leq k+10$,
then for any $u,\,v\in \bar{F}^{-3/4}$,

$$
\Big\|(i+\tau-p(\xi))^{-1}\eta_{k_1}(\xi)i\xi\widehat{P_{k}u}*\widehat{P_{k_2}v}
\Big\|_{X_{k_1}}\lesssim 2^{-3k/4}\|\widehat{P_{
k}u}\|_{X_{k}}\|\widehat{P_{k_2}v}\|_{X_{k_2}}.
$$
\end{lemma}

\begin{proof} As in the proof of Lemma \ref{h-l}, we may assume
$k=k_2$. Then
\begin{eqnarray}
&&\Big\|(i+\tau-p(\xi))^{-1}\eta_{k_1}(\xi)i\xi\widehat{P_{k}u}*\widehat{P_{k_2}v}
\Big\|_{X_{k_1}}\nonumber\\
&&\lesssim 2^{k_1}\sum_{j_1,j_2,j_3\geq
0}2^{-j_1/2}\|\chi_{D_{k_1,j_1}}\cdot u_{k,j_2}*v_{k,j_3}\|_{L^2_{\xi,\tau}},
\end{eqnarray}
where $u_{k,j_1},\,v_{k,j_2}$ are as in (\ref{u-v}) with $j_{max}\geq 3k-20$ and $j_i\leq 10k\,(i=1,2,3)$ being assumed in the summation. Applying (\ref{estimate1}) we
get
\begin{eqnarray*}
&&2^{k_1}\sum_{j_1,j_2,j_3\geq 0}2^{-j_1/2}\|\chi_{D_{k_1,j_1}}\cdot
u_{k,j_2}*v_{k,j_3}\|_{L^2_{\xi,\tau}}\\
&&\lesssim
\Big(\sum_{j_1=j_{max}}+\sum_{j_2=j_{max}}+\sum_{j_3=j_{max}}\Big)2^{-\frac{j_1}2}
2^{\frac{3k}4}2^{\frac{j_{min}}{2}}2^{\frac{j_{med}}4}\|u_{k,j_2}\|_{L^2_{\xi,\tau}}\|v_{k,j_3}\|_{L^2_{\xi,\tau}}\\
&&:=I+II+III.
\end{eqnarray*}
Since it is easy to get the bound of $I$ and there exists a symmetric relation between $II$ and $III$, it is enough to bound $II$ according to
$II\lesssim II_1+II_2$, where
\begin{eqnarray*}
&II_1=\sum_{j_2=j_{max},\,j_1\leq
j_3}2^{-\frac{j_1}2}2^{\frac{3k}4}2^{\frac{j_{min}}2}2^{\frac{j_{med}}4}\|u_{k,j_2}\|_{L^2_{\xi,\tau}}\|v_{k,j_3}\|_{L^2_{\xi,\tau}};\\
&II_2=\sum_{j_2=j_{max},\,j_1\geq
j_3} 2^{-\frac{j_1}2}2^{\frac{3k}4}2^{\frac{j_{min}}2}2^{\frac{j_{med}}4}\|u_{k,j_2}\|_{L^2_{\xi,\tau}}\|v_{k,j_3}\|_{L^2_{\xi,\tau}}.
\end{eqnarray*}
For $II_1$, by summing on $j_1$ we have
\begin{eqnarray*}
II_1&\lesssim& \sum_{j_2=j_{max},\,j_1\leq
j_3}2^{-j_1/2}2^{3k/4}2^{j_1/2}2^{j_3/4}\|u_{k,j_2}\|_{L^2_{\xi,\tau}}\|v_{k,j_3}\|_{L^2_{\xi,\tau}}\\
&\lesssim&\sum_{j_2\geq 3k-20,\,j_3\geq
0}2^{3k/4}2^{j_3/2}\|u_{k,j_2}\|_{L^2_{\xi,\tau}}\|v_{k,j_3}\|_{L^2_{\xi,\tau}}\\
&\lesssim&
2^{-3k/4}\|\widehat{P_ku}\|_{X_k}\|\widehat{P_{k_2}v}\|_{X_{k_2}}.
\end{eqnarray*}
For $II_2$, we have
\begin{eqnarray*}
II_2&\lesssim& \sum_{j_2=j_{max},\,j_1\geq
j_3}2^{-j_1/2}2^{3k/4}2^{j_3/2}2^{j_4/4}\|u_{k,j_2}\|_{L^2_{\xi,\tau}}\|v_{k,j_3}\|_{L^2_{\xi,\tau}}\\
&\lesssim&
2^{-3k/4}\|\widehat{P_ku}\|_{X_k}\|\widehat{P_{k_2}v}\|_{X_{k_2}}.
\end{eqnarray*}
\end{proof}

To consider the low-low interaction, from now on let $\psi\in
C_0^\infty(\mathbb R)$ be a standard bump function such that
$\psi(t)\equiv 1$ if $|t|<1$ and $\psi(t)\equiv 0$ if $|t|>2$.

\begin{lemma}[low -- low interaction]\label{l-l} If $0\leq k_1,\,k_2,\,k_3\leq 100$,
then for any $u,\,v\in F^s$,
$$\Big\|(i+\tau-p(\xi))^{-1}\eta_{k_1}(\xi)i\xi\widehat{\psi(t)P_{k_2}u}*\widehat{P_{k_3}v}
\Big\|_{X_{k_1}}\lesssim \|P_{
k_2}u\|_{L^\infty_tL^2_x}\|P_{k_3}v\|_{L^\infty_tL^2_x}.$$
\end{lemma}

\begin{proof} From the definition of $X_k$, Plancherel's
equality and Bernstein's inequality we achieve
\begin{eqnarray*}
&&\Big\|(i+\tau-p(\xi))^{-1}\eta_{k_1}(\xi)i\xi\widehat{\psi(t)P_{k_2}u}*\widehat{P_{k_3}v}
\Big\|_{X_{k_1}}\\
&&\lesssim 2^{k_1}\sum_{j_3\geq
0}2^{-j_3/2}\|\psi(t)P_{k_2}u*P_{k_3}v\|_{L^2_tL^2_x}\\
&&\lesssim
\|P_{k_2}u\|_{L^\infty_tL^2_x}\|P_{k_3}v\|_{L^\infty_tL^2_x},
\end{eqnarray*}
thereby reaching the desired estimate.
\end{proof}

\begin{lemma}[high -- high interaction]\label{h-h}
\item{\rm (i)} If $k\geq 10,\,
|k-k_2|\leq 5$, then for any $u,\,v\in F^s$,
\begin{equation}\label{hh1}\Big\|(i+\tau-p(\xi))^{-1}\eta_0(\xi)i\xi\widehat{P_{k}u}*\widehat{P_{k_2}v}
\Big\|_{X_0}\lesssim k2^{-3k/2}\|P_{
k}u\|_{X_k}\|P_{k_2}v\|_{X_{k_2}}.
\end{equation}

\item{\rm(ii)} If $k\geq 10,\, |k-k_2|\leq 5$ and $1\leq k_1\leq k-9$, then
for any $u,\,v\in F^s$,

\begin{eqnarray}\label{hh2}
&&\Big\|(i+\tau-p(\xi))^{-1}\eta_{k_1}(\xi)i\xi\widehat{P_{k}u}*\widehat{P_{k_2}v}
\Big\|_{X_1}\nonumber\\
&&\lesssim (2^{-3k/2}+k2^{-2k+k_1/2})\|P_{
k}u\|_{X_k}\|P_{k_2}v\|_{X_{k_2}}.
\end{eqnarray}
\end{lemma}
\begin{proof} For part (i), we may assume $k=k_2$. The
left-hand side of (\ref{hh1}) is dominated by
$$
\sum_{k_3=-\infty}^02^{k_3}\sum_{j_1,\,j_2,\,j_3\geq
0}2^{-j_3/2}\|\chi_{D_{k_3,\,j_3}}\cdot
u_{k,\,j_1}*v_{k,\,j_2}\|_{L^2_{\xi,\tau}},
$$
where $u_{k,\,j_1},\,v_{k,\,j_2}$ are as in (\ref{u-v}) with
$k_3\geq -10k$ and $j_1,\,j_2,\,j_3\leq 10k$ being assumed in the summation. Now, it
suffices to consider the worst case $|j_3-2k-k_3|\leq 10$: Applying (\ref{estimate2}) we get
\begin{eqnarray*}
&&\Big\|(i+\tau-p(\xi))^{-1}\eta_0(\xi)i\xi\widehat{P_{k}u}*\widehat{P_{k_2}v}
\Big\|_{X_0}\\
&&\lesssim\sum_{k_3=-10k}^0\sum_{j_1,\,j_2\geq
0}2^{-k}2^{-k_3/2}2^{k_3}2^{-k/2}2^{-k_3/2}2^{j_1/2}2^{j_2/2}\|u_{k,\,j_1}\|_{L^2_{\xi,\tau}}
\|v_{k,\,j_2}\|_{L^2_{\xi,\tau}}\\
&&\lesssim
k2^{-3k/2}\|\widehat{P_{k}u}\|_{X_k}\|\widehat{P_{k}v}\|_{X_k},
\end{eqnarray*}
which is the desired estimate in part (a).

For part (ii) we may also assume $k=k_2$, and consequently get
\begin{eqnarray}\label{hh3}
&&\Big\|(i+\tau-p(\xi))^{-1}\eta_{k_1}(\xi)i\xi\widehat{P_{k}u}*\widehat{P_{k_2}v}
\Big\|_{X_1}\nonumber\\
&&\lesssim 2^{k_1}\sum_{j_1,\,j_2,\,j_3\geq
0}2^{-j_1/2}\|\chi_{D_{k_1,\,j_1}}\cdot
u_{k,\,j_2}*v_{k,\,j_3}\|_{L^2_{\xi,\tau}},
\end{eqnarray}
where $u_{k,\,j_2},\,v_{k,\,j_3}$ are as in (\ref{u-v}) with
$k_{max}\geq 2k+k_1-10$ and $j_1,\,j_2,\,j_3\leq 10k$ being assumed in the summation.
We will bound the right-hand side of (\ref{hh3}) case by case. The
first case is that $j_1=j_{max}$ in the summation. Concerning this case we apply
(\ref{estimate2}) to get
\begin{eqnarray*}
&&2^{k_1}\sum_{j_1,\,j_2,\,j_3\geq
0}2^{-j_1/2}\|\chi_{D_{k_1,\,j_1}}\cdot u_{k,\,j_2}*v_{k,\,j_3}\|_{L^2_{\xi,\tau}}\\
&&\lesssim 2^{k_1}\sum_{j_1\geq 2k+k_1-10}\sum_{j_2,\,j_3\geq
0}2^{-j_1}2^{-k/2}2^{-k_1/2}2^{j_2/2}2^{j_3/2}\|u_{k,\,j_2}\|_{L^2_{\xi,\tau}}
\|v_{k,\,j_3}\|_{L^2_{\xi,\tau}}\\
&&\lesssim
2^{-3k/2}\|\widehat{P_{k}u}\|_{X_k}\|\widehat{P_{k}v}\|_{X_k},
\end{eqnarray*}
which is acceptable. The second case is $j_2=j_{max}$. Regarding this one we apply (\ref{estimate3}) to get
\begin{eqnarray*}
&&2^{k_1}\sum_{j_1,\,j_2,\,j_3\geq
0}2^{-j_1/2}\|\chi_{D_{k_1,\,j_1}}\cdot u_{k,\,j_2}*v_{k,\,j_3}\|_{L^2_{\xi,\tau}}\\
&&\lesssim 2^{k_1}\sum_{j_2\geq 2k+k_1-10}\sum_{j_1,\,j_3\geq
0}2^{-j_1}2^{-k}2^{-k_1/2}2^{j_1/2}2^{j_2/2}\|u_{k,\,j_2}\|_{L^2_{\xi,\tau}}
\|v_{k,\,j_3}\|_{L^2_{\xi,\tau}}\\
&&\lesssim
k2^{-2k}2^{k_1/2}\|\widehat{P_{k}u}\|_{X_k}\|\widehat{P_{k_2}v}\|_{X_{k_2}},
\end{eqnarray*}
where in the last inequality we have used $j_1\leq 10k$. The third case is
$j_3=j_{max}$, but this is identical with the case $j_2=j_{max}$ due to symmetry.
Thus, the estimate of (\ref{hh2}) is done.
\end{proof}

The main reason of using $\bar{F}^{-3/4}$ is the logarithmic
loss of derivative in (\ref{hh1}). Nevertheless, we can avoid the
logarithmic loss in (\ref{hh1}) by using a $\bar{X}_0$ structure on
the low frequency.

\begin{lemma}[$\bar{X}_0$ estimate]\label{propXbares}
Let $|k_1-k_2|\leq 5$ and $k_1\geq 10$. Then we have for all $u,v\in
\bar{F}^0$
\begin{align*}
\normo{\psi(t)\int_0^t W(t-s)P_{\leq
0}\partial_x[P_{k_1}u(s)P_{k_2}v(s)]ds}_{L_x^2 L_t^\infty}\les
2^{-\half{3k_1}}\norm{\wh{P_{k_1}u}}_{X_{k_1}}\norm{\wh{P_{k_2}u}}_{X_{k_2}}.
\end{align*}
\end{lemma}

\begin{proof} From now on we put
$$
Q(u,v)=\psi(t)\int_0^t W(t-s)P_{\leq
0}\partial_x[P_{k_1}u(s)P_{k_2}v(s)]ds.
$$
Via a straightforward computation we find a constant $c$ such that
\begin{align*}
\ft\left[Q(u,v)\right](\xi,\tau)=&c\int_\R
\frac{\widehat{\psi}(\tau-\tau')-\widehat{\psi}(\tau-p(\xi))}{\tau'-p(\xi)}\eta_0(\xi)i\xi
\\
&\times\
d\tau'\int_{\xi=\xi_1+\xi_2,\tau'=\tau_1+\tau_2}\wh{P_{k_1}u}(\xi_1,\tau_1)\wh{P_{k_2}v}(\xi_2,\tau_2).
\end{align*}
For the fixed point $\xi \in \R$ we split the hyperplane
$$
\Gamma:=\{\xi=\xi_1+\xi_2,\tau'=\tau_1+\tau_2\}
$$
into
\begin{align*}
\Gamma_1=&\ \{ |\xi|\les 2^{-2k_1}\} \cap \Gamma;\\
\Gamma_2=&\ \{ |\xi|\gg 2^{-2k_1},
|\tau_i-p(\xi_i)|\ll 3\cdot 2^{2k_1}|\xi|, i=1,2\} \cap \Gamma;\\
\Gamma_3=&\ \{ |\xi|\gg 2^{-2k_1},
|\tau_1-p(\xi_1)|\ges 3\cdot 2^{2k_1}|\xi|\} \cap \Gamma;\\
\Gamma_4=&\ \{ |\xi|\gg 2^{-2k_1}, |\tau_2-p(\xi_2)|\ges 3\cdot
2^{2k_1}|\xi|\} \cap \Gamma.
\end{align*}
With this splitting we write
\[\ft\left[\psi(t)\cdot
\int_0^tW(t-s)P_{\leq
0}\partial_x[P_{k_1}u(s)P_{k_2}v(s)]ds\right](\xi,\tau):=G_1+G_2+G_3+G_4,
\]
where
\[
G_i=C\int_\R
\frac{\widehat{\psi}(\tau-\tau')-\widehat{\psi}(\tau-p(\xi))}{\tau'-p(\xi)}\eta_0(\xi)i\xi
\int_{\Gamma_i}\wh{P_{k_1}u}(\xi_1,\tau_1)\wh{P_{k_2}v}(\xi_2,\tau_2)d\tau'
\]
with $C$ being a constant.

First of all, let us deal with the contribution of $G_1$. Using Lemma
\ref{lem2} and Proposition \ref{proplineares} (ii) (in the coming next section), we get
\[
\norm{\ft^{-1}(G_1)}_{ L_x^2L_t^\infty}\lesssim \normo{(i+\tau'-p(\xi))^{-1}\eta_0(\xi)i\xi
\int_{\Gamma_1}\wh{P_{k_1}u}(\xi_1,\tau_1)\wh{P_{k_2}v}(\xi_2,\tau_2)}_{X_0}.
\]
Since $|\xi|\lesssim 2^{-2k_1}$ holds in the region of $G_1$, we get
\begin{align*}
&\normo{(i+\tau'-p(\xi))^{-1}\eta_0(\xi)i\xi
\int_{G_1}\wh{P_{k_1}u}(\xi_1,\tau_1)\wh{P_{k_2}v}(\xi_2,\tau_2)}_{X_0}\\
&\lesssim \sum_{k_3\leq -2k_1+10}\sum_{j_3\geq 0}
2^{-j_3/2}2^{k_3}\sum_{j_1\geq 0, j_2\geq
0}\norm{\chi_{D_{k_3,j_3}}\cdot u_{k_1,j_1}*v_{k_2,j_2}}_{L^2}
\end{align*}
where
\[
u_{k_1,j_1}(\xi,\tau)=\eta_{k_1}(\xi)\eta_{j_1}(\tau-p(\xi))\wh{u}(\xi,\tau),
v_{k_1,j_1}(\xi,\tau)=\eta_{k_1}(\xi)\eta_{j_1}(\tau-p(\xi))\wh{v}(\xi,\tau).
\]
Using Lemma \ref{block estimate} (iii) we obtain
\begin{align*}
\norm{\ft^{-1}(G_1)}_{ L_x^2L_t^\infty}\lesssim& \sum_{k_3\leq
-2k_1+10}\sum_{j_i\geq 0} 2^{-j_3/2}2^{k_3}2^{j_{min}/2}2^{k_3/2}
\norm{u_{k_1,j_1}}_{L^2}\norm{v_{k_2,j_2}}_{L^2}\\
\lesssim&\ 2^{-3k_1}\norm{\wh{P_{k_1}u}}_{X_{k_1}}\norm{\wh{P_{k_2}u}}_{X_{k_2}},
\end{align*}
which suffices to give the bound for $G_1$.

Next, we settle the contribution of $G_3$. Using Lemma \ref{lem2} and Proposition \ref{proplineares}(ii) we get
\begin{align*}
\norm{\ft^{-1}(G_3)}_{ L_x^2L_t^\infty}\lesssim&\ \normo{(i+\tau'-p(\xi))^{-1}\eta_0(\xi)i\xi \int_{\Gamma_3}\wh{P_{k_1}u}(\xi_1,\tau_1)\wh{P_{k_2}v}(\xi_2,\tau_2)}_{X_0}\\
\lesssim&\ \sum_{k_3\leq 0}\sum_{j_3\geq 0} 2^{-j_3/2}2^{k_3}\sum_{j_1\geq
0, j_2\geq 0}\norm{\chi_{D_{k_3,j_3}}\cdot
u_{k_1,j_1}*v_{k_2,j_2}}_{L^2}
\end{align*}
Clearly we may assume $j_3\leq 10k_1$ in the summation above. Applying
Lemma 2.3 (iii) we get
\begin{align*}
\norm{\ft^{-1}(G_3)}_{ L_x^2L_t^\infty}\lesssim& \sum_{k_3\leq
0}\sum_{j_1\geq k_3+2k_1-10, j_2,j_3\geq 0}
2^{k_3}2^{j_{2}/2}2^{-3k_1}
\norm{u_{k_1,j_1}}_{L^2}\norm{v_{k_2,j_2}}_{L^2}\\
\lesssim&\ k_12^{-2k_1}\norm{\wh{P_{k_1}u}}_{X_{k_1}}\norm{\wh{P_{k_2}u}}_{X_{k_2}},
\end{align*}
which suffices to give the bound for $G_3$. From symmetry,
the bound for $G_4$ is the same as $G_3$.

Finally, we consider the contribution of $G_2$. From the proof
of the dyadic bilinear estimates, we know this term is the main
contribution. By a computation we get
\begin{align*}
\ft_t^{-1}(G_2)=\psi(t)\int_0^t \frac{\eta_0(\xi)i\xi}{e^{i(s-t)p(\xi)}}
\int_{\R^2}e^{is(\tau_1+\tau_2)}
\int_{\xi=\xi_1+\xi_2}\frac{{u_{k_1}}(\xi_1,\tau_1)}{[{v_{k_2}}(\xi_2,\tau_2)]^{-1}}\
d\tau_1 d\tau_2ds
\end{align*}
where
\begin{align*}
u_{k_1}(\xi_1,\tau_1)=&\ \eta_{k_1}(\xi_1)\chi_{\{|\tau_1-p(\xi_1)|\ll
3\cdot 2^{2k_1}|\xi|\}}\wh{u}(\xi_1,\tau_1),\\
v_{k_2}(\xi_2,\tau_2)=&\ \eta_{k_2}(\xi_2)\chi_{\{|\tau_2-p(\xi_2)|\ll
3\cdot 2^{2k_1}|\xi|\}}\wh{v}(\xi_2,\tau_2).
\end{align*}
By a change of variable $\tau_1'=\tau_1-p(\xi_1)$,
$\tau_2'=\tau_2-p(\xi_2)$, we get
\begin{align*}
\ft_t^{-1}(G_2)=&\ \psi(t)e^{itp(\xi)}\eta_0(\xi)\xi\int_{\R^2}e^{it(\tau_1+\tau_2)}
\int_{\xi=\xi_1+\xi_2}\frac{e^{it(p(\xi_1)+p(\xi_2)-p(\xi))}-e^{-it(\tau_1+\tau_2)}}
{\tau_1+\tau_2-p(\xi)+p(\xi_1)+p(\xi_2)}\\
&\times \
{u_{k_1}}(\xi_1,\tau_1+p(\xi_1)){v_{k_2}}(\xi_2,\tau_2+p(\xi_2))\
d\tau_1 d\tau_2\\
:=&\ \ft_t^{-1}(I)-\ft_t^{-1}(II).
\end{align*}
For the contribution of $\ft_t^{-1}(II)$, we have
\begin{align*}
&\ft_t^{-1}(II)=\int_{\R^2}\frac{\psi(t)\eta_0(\xi)\xi}{e^{-itp(\xi)}}\int_{\xi=\xi_1+\xi_2}
\frac{{u_{k_1}}(\xi_1,\tau_1+p(\xi_1)){v_{k_2}}(\xi_2,\tau_2+p(\xi_2))}{\tau_1+\tau_2-p(\xi)+p(\xi_1)+p(\xi_2)}\
d\tau_1 d\tau_2.
\end{align*}
Since in the support of $u_{k_1}$ and $u_{k_2}$ we have
$$
|\tau_1+\tau_2-p(\xi)+p(\xi_1)+p(\xi_2)|\sim 2^{2k_1}|\xi|,
$$
we conclude from Lemma \ref{lem1} that
\begin{align*}
\norm{\ft^{-1}(II)}_{L_x^2L_t^\infty}\lesssim&\ \int_{\R^2}\normo{\int_{\xi=\xi_1+\xi_2}\xi\frac{{u_{k_1}}(\xi_1,\tau_1+p(\xi_1)){v_{k_2}}(\xi_2,\tau_2+p(\xi_2))}
{\tau_1+\tau_2-p(\xi)+p(\xi_1)+p(\xi_2)}}_{L_\xi^2}d\tau_1d\tau_2\\
\lesssim&\ 2^{-\half{3k_1}}\norm{\wh{P_{k_1}u}}_{X_{k_1}}\norm{\wh{P_{k_2}u}}_{X_{k_2}}.
\end{align*}

To complete the argument, it remains to prove the following inequality
\begin{eqnarray*}
\norm{\ft^{-1}(I)}_{L_x^2L_t^\infty}\lesssim
2^{-3k_1/2}\norm{\wh{P_{k_1}u}}_{X_{k_1}}\norm{\wh{P_{k_2}u}}_{X_{k_2}}.
\end{eqnarray*}
In doing this, let us compare the term $I$ with the following term $I'$:
\begin{align*}
\ft_t^{-1}(I')
=&\ \psi(t)e^{itp(\xi)}\eta_0(\xi)\xi\int_{\R^2}e^{it(\tau_1+\tau_2)}\int_{\xi=\xi_1+\xi_2}
\frac{e^{it(p(\xi_1)+p(\xi_2)-p(\xi))}}{-p(\xi)+p(\xi_1)+p(\xi_2)}\\
&\times \
{u_{k_1}}(\xi_1,\tau_1+p(\xi_1)){v_{k_2}}(\xi_2,\tau_2+p(\xi_2))\
d\tau_1 d\tau_2.
\end{align*}
By symmetry we may assume that $k_1\geq k_2$ as well as $|\alpha|\leq
2^{k_1}|\beta|$. Then, on the hyperplane $\xi=\xi_1+\xi_2$ we have
\[-p(\xi_1+\xi_2)+p(\xi_1)+p(\xi_2)=3\beta\xi_1\xi_2\xi\big(1-\frac{2\alpha}{3\beta|\xi_1|}\big),\]
whence getting
\[\frac{1}{1-\frac{2\alpha}{3\beta|\xi_1|}}=\sum_{n=0}^\infty\brk{\frac{2\alpha}{3\beta|\xi_1|}}^n.\]
Inserting this into $I'$ we have
\begin{align*}
\ft_t^{-1}(I')
=&\ \psi(t)\eta_0(\xi)\sum_{n=0}^\infty\int_{\R^2}e^{it(\tau_1+\tau_2)}\int_{\xi=\xi_1+\xi_2}
e^{it(p(\xi_1)+p(\xi_2))}\frac{(2\alpha/3\beta)^n}{(\xi_1\xi_2)|\xi_1|^{n}}\\
&\times \
{u_{k_1}}(\xi_1,\tau_1+p(\xi_1)){v_{k_2}}(\xi_2,\tau_2+p(\xi_2))\
d\tau_1 d\tau_2.
\end{align*}
Since it is easy to see that (actually we need a smooth version of
$\chi_{\{|\xi|\gg \lambda\}}$): Any $\lambda>0$ ensures 
\[\norm{\ft_x^{-1}\chi_{\{|\xi|\gg \lambda\}}\ft_x u}_{L_x^2L_t^\infty}\lesssim \norm{u}_{L_x^2L_t^\infty},
\]
by setting
\[\ft(f_{\tau_1})(\xi)=\wh{P_{k_1}u}(\xi,\tau_1+p(\xi));\ \ \ft(g_{\tau_2})(\xi)=\wh{P_{k_2}v}(\xi,\tau_2+p(\xi)),\] 
we get from Lemma \ref{lem1} that 
\begin{align*}
\norm{\ft^{-1}(I')}_{L_x^2L_t^\infty}\lesssim&\ \sum_{n=0}^\infty C^n
\int_{\R^2}\norm{W(t)\partial_x^{-{(n+1)}}f_{\tau_1}W(t)\partial_x^{-{1}}g_{\tau_2}}_{L_x^2L_t^\infty}d\tau_1d\tau_2\\
\lesssim&\ \sum_{n=0}^\infty C^n\int_{\R^2}\norm{W(t)\partial_x^{-{(n+1)}}f_{\tau_1}}_{L_x^4L_t^\infty}\norm{W(t)\partial_x^{-{1}}g_{\tau_2}}_{L_x^4L_t^\infty}d\tau_1d\tau_2\\
\lesssim&\ 2^{-\half{3k_1}}\norm{\wh{P_{k_1}u}}_{X_{k_1}}\norm{\wh{P_{k_2}u}}_{X_{k_2}}.
\end{align*}
Meanwhile, it is also necessary to establish the following inequality
\begin{align*}
\norm{\ft^{-1}(I-I')}_{L_x^2L_t^\infty}\lesssim
2^{-\half{3k_1}}\norm{\wh{P_{k_1}u}}_{X_{k_1}}\norm{\wh{P_{k_2}u}}_{X_{k_2}}.
\end{align*}
Since in the integral region we have $|\tau_i|\ll 2^{2k_1}|\xi|$ where
$i=1,2$, on the hyperplane $\xi=\xi_1+\xi_2$ we obtain
\begin{align*}
&\frac{1}{\tau_1+\tau_2-p(\xi)+p(\xi_1)+p(\xi_2)}-\frac{1}{-p(\xi)+p(\xi_1)+p(\xi_2)}\\
&=\ \sum_{n=1}^\infty
\rev{-p(\xi)+p(\xi_1)+p(\xi_2)}\brk{\frac{\tau_1+\tau_2}{-p(\xi)+p(\xi_1)+p(\xi_2)}}^n\\
&=\ C\sum_{n=1}^\infty\frac{1}{\xi_1\xi_2\xi}\sum_{k=0}^\infty
\brk{\frac{2\alpha}{3\beta|\xi_1|}}^k
\brk{\frac{\tau_1+\tau_2}{\xi_1\xi_2\xi}}^n \sum_{j_1,\cdots,
j_n=0}^\infty
\prod_{i=1}^n\brk{\frac{2\alpha}{3\beta|\xi_1|}}^{j_i}.
\end{align*}
The purpose of doing such a decomposition is to make the variable separately. So, we can apply Lemma \ref{lem1}. Via decomposing low frequency we get
\begin{align*}
&\ft_t^{-1}(I-I')=\sum_{n=1}^\infty\psi(t)\eta_0(\xi)\int_{\R^2}e^{it(\tau_1+\tau_2)}\sum_{2^{k_3}\gg
2^{-2k_1}\max(|\tau_1|,|\tau_2|)}\chi_{k_3}(\xi)\\
&\times\int_{\xi=\xi_1+\xi_2}e^{it(p(\xi_1)+p(\xi_2))}
{u_{k_1}}(\xi_1,\tau_1+p(\xi_1){v_{k_2}}(\xi_2,\tau_2+p(\xi_2)\frac{1}{\xi_1\xi_2}\\
&\times \sum_{k=0}^\infty
\brk{\frac{2\alpha}{3\beta|\xi_1|}}^k
\brk{\frac{\tau_1+\tau_2}{\xi_1\xi_2\xi}}^n \sum_{j_1,\cdots,
j_n=0}^\infty
\prod_{i=1}^n\brk{\frac{2\alpha}{3\beta|\xi_1|}}^{j_i}
d\tau_1 d\tau_2.
\end{align*}
Using the fact that $\chi_{k_3}(\xi)(\xi/2^{k_3})^{-n}$ is a
multiplier for the space $L_x^2L_t^\infty$, we get
\begin{align*}
&\norm{\ft^{-1}(I-I')}_{L_x^2L_t^\infty}\\
&\lesssim\ \sum_{n=1}^\infty \int_{\R^2} \sum_{2^{k_3}\gg
2^{-2k_1}\max(|\tau_1|,|\tau_2|)} C^n
|\tau_1+\tau_2|^n2^{-nk_3}2^{-2nk_1}\\
&\ \times \
2^{-3k_1/2}\norm{\ft(f_{\tau_1})}_{L^2}\norm{\ft(g_{\tau_2})}_{L^2}d\tau_1
d\tau_2\\
&\lesssim\ 2^{-3k_1/2}\norm{\wh{P_{k_1}u}}_{X_{k_1}}\norm{\wh{P_{k_2}u}}_{X_{k_2}}.
\end{align*}
\end{proof}

\section{Proof of Theorem \ref{thm3}}

This section is devoted to proving Theorem \ref{thm3} by using
the standard contraction principle and the dyadic estimates obtained in the last section.

While making a comparison with the KdV equation, we will immediately encounter an essential difference -- unlike KdV, the Benjamin equation has no any scaling invariant property. Instead, in order to prove Theorem \ref{thm3} we
may use the following scale argument: If $u(x,t)$ is a solution of
(\ref{Ben1}), then for $\lambda>0$,
$$
u_\lambda(x,t)=\lambda^2 u(\lambda x, \lambda^3t)
$$
is a solution to the following equation
\begin{align}\label{sBen}
\left\{
\begin{array}{l}
\partial_t u - \tilde\gamma\partial_x u+\tilde\alpha {\mathcal
H}\partial^2_xu+\beta\partial^3_xu+\partial_x(u^2) =0, \quad
(x,t)\in {\mathbb R\times\mathbb
R},\\
u(x,0) =\phi(x),\quad x\in\mathbb R,
\end{array}
\right.
\end{align}
where
$$
\tilde\gamma=\lambda^2\gamma,\quad\tilde\alpha=\lambda\alpha,\quad\hbox{
and}\quad \phi(x)=\lambda^2u_0(\lambda x).
$$

Note that
$$
\norm{\lambda^2u_0(\lambda
x)}_{{H}^s}\leq\lambda^{s+3/2}\norm{u_0}_{{H}^s}+\lambda^{3/2}\norm{u_0}_{{H}^s}.
$$
So, under $s>-3/2$ we may assume $\norm{\phi}_{H^s}\ll 1$ by taking
$0<\lambda\ll 1$. Also, upon a normalization of $u$, we may assume
$\beta=1$. With these assumptions, we see that a consideration of the local well-posedness of (\ref{Ben1}) in
$H^{-3/4}(\mathbb R)$ is reduced to handling the similar matter for the equation (\ref{sBen}) under
the condition
$$
|\tilde\gamma|\leq 1,\,|\tilde\alpha|\leq 1, \quad
\norm{\phi}_{H^{-3/4}}\ll 1.
$$

Using the Duhamel principle and setting $u(t)=u(x,t)$ we see that (\ref{sBen}) can be solved by finding the unique solution of the following truncated integral equation
$$
u(t)=\psi\big(\frac{t}{4}\big)\Big[W(t)\phi-\int_0^t W(t-\tau)\partial_x(\psi^2(\tau)u^2(\tau))\,d\tau\Big].
$$
In solving this last equation, we need the forthcoming ingredients.

\begin{proposition}[linear estimates]\label{proplineares}
\item{\rm (i)} If $s\in \R$ and  $\phi \in H^s$, then there exists $C>0$
such that
\begin{eqnarray}
\norm{\psi(t)W(t)\phi}_{\bar{F}^s}\leq C\norm{\phi}_{H^{s}}.
\end{eqnarray}

\item{\rm(ii)} If $s\in \R, k\in \Z_+$ and $u$ satisfies
$(i+\tau-p(\xi))^{-1}\ft(u)\in X_k$, then there exists $C>0$ such
that
\begin{eqnarray}
\normo{\ft\left[\psi(t)\int_0^t W(t-s)(u(s))ds\right]}_{X_k}\leq
C\norm{(i+\tau-p(\xi))^{-1}\ft(u)}_{X_k}.
\end{eqnarray}
\end{proposition}

\begin{proof} A proof of (i) follows from Lemma \ref{lem1}. A proof of (ii) can be given via \cite{GuoWang}.
\end{proof}

\begin{proposition}[bilinear estimates]\label{propbilinearbd} For $u,v\in \bar{F}^s$ let
\begin{align}
B(u,v):=\psi\big(\frac{t}{4}\big)\int_0^tW(t-\tau)\partial_x\big(\psi^2(\tau)u(\tau)\cdot
v(\tau)\big)d\tau.
\end{align}
If $-3/4\leq s\leq 0$, then there exists $C>0$ such that
\begin{eqnarray}\label{eq:bilinearbd}
\norm{B(u,v)}_{\bar{F}^s}\leq
C(\norm{u}_{\bar{F}^s}\norm{v}_{\bar{F}^{-3/4}}+\norm{u}_{\bar{F}^{-3/4}}\norm{v}_{\bar{F}^s})
\end{eqnarray}
hold for any $u,v\in \bar{F}^s$.
\end{proposition}
\begin{proof} In light of the argument for \cite[Proposition 4.2]{Guo}, we check the proposition as follows.
Thanks to
\begin{equation}\label{eqA}
\|B(u,v)\|^2=\|P_{\le 0}B(u,v)\|^2_{\bar{X}_0}+\sum_{k_1\ge 1}2^{2k_1s}\|\eta_{k_1}(\xi)\mathcal{F}[B(u,v)]\|^2_{X_{k_1}},
\end{equation}
we are about to control the two terms of the right-hand side of (\ref{eqA}).

Using the decomposition of $u,v$ we have
$$
\|B(u,v)\|_{\bar{X}_0}\le\sum_{k_2,k_3\ge 0}\|B(P_{k_2}u,P_{k_3}v)\|_{\bar{X}_0},
$$
thereby considering two cases:

(i) If $\max(k_2,k_3)\le 10$, then Lemma \ref{lem2} implies
$$
\|\eta_0(t)P_{\le 0}u\|_{\bar{X}_0}\lesssim \|P_{\le 0}u\|_{\bar{X}_0}.
$$
This, along with Lemma \ref{l-l} and Proposition \ref{proplineares}, gives
$$
\|B(P_{k_2}u,P_{k_3}v)\|_{\bar{X}_0}\lesssim \|P_{k_2}u\|_{L^\infty_tL_x^2}
\|P_{k_3}v\|_{L^\infty_tL_x^2},
$$
whence yielding
\begin{equation}\label{eqB}
\|B(u,v)\|_{\bar{X}_0}\lesssim(\norm{u}_{\bar{F}^s}\norm{v}_{\bar{F}^{-3/4}}+\norm{u}_{\bar{F}^{-3/4}}\norm{v}_{\bar{F}^s}).
\end{equation}

(ii) If $\max(k_2,k_3)>10$, then $|k_2-k_3|\le 5$ and hence by Lemma
\ref{propXbares},
\begin{eqnarray}\label{eqC}
\|B(u,v)\|_{\bar{X}_0}&\le&\sum_{|k_2-k_3|\le 5,k_2,k_3\ge 10}2^{-3k_2/2}\|\mathcal F(P_{k_2}u)\|_{X_{k_2}}\|\mathcal F(P_{k_3}v)\|_{X_{k_3}}\nonumber\\
&\lesssim&\norm{u}_{\bar{F}^{-3/4}}\norm{v}_{\bar{F}^{-3/4}}\\
&\lesssim&(\norm{u}_{\bar{F}^s}\norm{v}_{\bar{F}^{-3/4}}+\norm{u}_{\bar{F}^{-3/4}}\norm{v}_{\bar{F}^s}).\nonumber
\end{eqnarray}

Now a combination of (\ref{eqB}) and (\ref{eqC}) deduces
\begin{equation}\label{eqT1}
\|P_{\le 0}B(u,v)\|_{\bar{X}_0}\lesssim\norm{u}_{\bar{F}^s}\norm{v}_{\bar{F}^{-3/4}}+\norm{u}_{\bar{F}^{-3/4}}\norm{v}_{\bar{F}^s}.
\end{equation}

Next, let us control the second part at the right-hand side of (\ref{eqA}). To do so, owing to symmetry we may assume $k_2\le k_3$. Decomposing $u$ and $v$ again and using Proposition \ref{proplineares} (ii), we see
\begin{eqnarray*}
&&\|\eta_{k_1}(\xi)\mathcal{F}[B(u,v)]\|_{X_{k_1}}\\
&&\lesssim\sum_{k_2,k_3\ge 0}\|\eta_{k_1}(\xi)\mathcal{F}[B(P_{k_2}u,P_{k_3}v)]\|_{X_{k_1}}\\
&&\lesssim\sum_{k_2,k_3\ge 0}\|(i+\tau-p(\xi))^{-1}\eta_{k_1}(\xi)i\xi\widehat{\psi(t)P_{k_2}u}*\widehat{\psi(t)P_{k_3}v}
\Big\|_{X_{k_1}}.
\end{eqnarray*}

(iii) If $k_{\max}\le 20$, then an application of Lemma \ref{l-l} and (\ref{eq05}) derives
\begin{eqnarray*}
&&\sum_{k_2, k_3\ge 0}\|(i+\tau-p(\xi))^{-1}\eta_{k_1}(\xi)i\xi\widehat{\psi(t)P_{k_2}u}*\widehat{\psi(t)P_{k_3}v}
\Big\|_{X_{k_1}}\\
&&\lesssim \sum_{k_{\max}\le 20}\|P_{k_2}u\|_{L_t^\infty L_x^2}\|P_{k_3}v\|_{L_t^\infty L_x^2}.
\end{eqnarray*}
Note that
\[
\|P_{k}v\|_{L_t^\infty L_x^2}\lesssim
\left\{
\begin{array}{l}

\|P_{k_3}v\|_{X_k}\quad\hbox{when}\quad k\ge 1,\\
\|P_{k_3}v\|_{\bar{X}_k}\quad\hbox{when}\quad k=0.
\end{array}
\right.
\]
So we get
\begin{eqnarray}\label{eqD}
&&\sum_{k_1\ge 1}2^{2k_1s}\Big[\sum_{k_2,k_3\ge 0}\|(i+\tau-p(\xi)){-1}\eta_{k_1}(\xi)i\xi\widehat{\psi(t)P_{k_2}u}*\widehat{\psi(t)P_{k_3}v}
\Big\|_{X_{k_1}}\Big]\nonumber\\
&&\lesssim (\norm{u}_{\bar{F}^{-3/4}}\norm{v}_{\bar{F}^s})^2.
\end{eqnarray}

(iv) If $k_{\max}>20$, then three subcases are considered:
\[
\left\{
\begin{array}{l}
\rm{(iv)_1}: |k_1-k_3|\le 5,\ k_2\le k_1-10;\\
\rm{(iv)_2}: |k_1-k_3|\le 5,\ k_1-9\le k_2\le k_3;\\
\rm{(iv)_3}: |k_2-k_3|\le 5,\ 1\le k_1\le k_2-5.
\end{array}
\right.
\]
For $\rm{(iv)_1}$, we use Lemma \ref{h-l} (i) with $k_2=0$ and Lemma \ref{h-l} (ii) with $k_2\ge 1$ to get (\ref{eqD}).
For $\rm{(iv)_2}$, we use Lemma \ref{l-l} to establish (\ref{eqD}).  For $\rm{(iv)_3}$, we apply Lemma \ref{h-h} (ii) to achieve (\ref{eqD}).

A combination of (iii) and (iv) implies
\begin{equation}\label{eqT2}
\sum_{k_1\ge 1}2^{2k_1s}\|\eta_{k_1}(\xi)\mathcal{F}[B(u,v)]\|^2_{X_{k_1}}\lesssim
\norm{u}_{\bar{F}^{-3/4}}\norm{v}_{\bar{F}^s}.
\end{equation}

Finally, we bring (\ref{eqT1}) and (\ref{eqT2}) into (\ref{eqA}) to produce the bilinear estimate (\ref{eq:bilinearbd}).
\end{proof}

Keeping the previous two propositions in mind, we can use the
standard fixed point argument (for the bounded bilinear operator
$B:\bar{F}^s\times\bar{F}^s\mapsto \bar{F}^s$ whenever $s\in [-3/4,0]$) to find a solution $u$ of (\ref{sBen}) in both $\bar{F}^{-3/4}$ and $C([-T,T]; H^{-3/4})$ for some $T>0$ depending on
the initial data $\phi$, and then verify the rest of Theorem \ref{thm3}.

\section{Modified Energies for Global Well-posedness}

To apply the so-called I-method \cite{I-method} to extending the
local solution to the global, let us review a couple of definitions. Given a complex-valued function $m: \R^k \rightarrow \C$, we say that $m$ is symmetric provided $m(\xi_1,\cdots, \xi_k)=m(\sigma(\xi_1,\cdots,
\xi_k))$ holds for all $\sigma \in S_k$, the group of all permutations on
$k$ objects. The symmetrization of $m$ is the function
\begin{equation*}
[m]_{sym}(\xi_1,\xi_2,\cdots, \xi_k)=\rev{k!}\sum_{\sigma\in
S_k}m(\sigma(\xi_1,\xi_2,\cdots,\xi_k)).
\end{equation*}
We then define a $k$-linear functional associated to the function $m$
(multiplier) acting on $k$ functions $u_1,\cdots,u_k$,
\begin{equation*}
\Lambda_k(m;u_1,\cdots,u_k)=\int_{\xi_1+\cdots+\xi_k=0}m(\xi_1,\cdots,\xi_k)\widehat{u_1}(\xi_1)\cdots
\widehat{u_k}(\xi_k).
\end{equation*}
In the sequel, we will often apply $\Lambda_k$ to $k$ copies of the same function
$u$. Consequently, $\Lambda_k(m;u,\ldots,u)$ may simply be written $\Lambda_k(m)$.
Using the symmetry of the measure on hyperplane, we obtain
$\Lambda_k(m)=\Lambda_k([m]_{sym})$, thereby reaching the following assertion.

\begin{lemma}[ODE in time]\label{lem51}
Suppose $u$ satisfies the Benjamin equation \eqref{sBen} and that
$m$ is a symmetric function. Then
\begin{align}\label{eq:menergy}
\frac{d}{dt}\Lambda_k(m)&=\Lambda_k(mv_k)-\Lambda_k(mh_k)\\
&-i\half k\Lambda_{k+1}(m(\xi_1,\ldots,\xi_{k-1},\xi_k+\xi_{k+1})(\xi_k+\xi_{k+1})),\nonumber
\end{align}
where
\[v_k=i(\xi_1^3+\xi_2^3+\cdots+\xi_k^3)\quad\hbox{and}\quad
h_k=i\tilde\alpha(\xi_1|\xi_1|+\xi_2|\xi_2|+\cdots+\xi_k|\xi_k|).\]
\end{lemma}
\begin{proof} This may be directly verified by the Benjamin equation \eqref{sBen}.
\end{proof}

Next, we define a branch of the modified energies. Given an arbitrary even
$\R$-valued function $m:\R\rightarrow \R$, let
\begin{align*}
\widehat{If}(\xi)=m(\xi)\widehat{f}(\xi),
\end{align*}
where the multiplier $m(\xi)$ is smooth, monotone, and of the form:
\begin{align}\label{eq:m}
m(\xi)=\left\{
\begin{array}{r}
1, \quad \quad |\xi|<N,\\
N^{-s}|\xi|^s,\quad  |\xi|>2N,
\end{array}
\right.
\end{align}
when $N\gg 1$. The modified energy $E_I^2(t)$ is determined by
\begin{align*}
E_I^2(t)=\norm{Iu(t)}_{L^2}^2.
\end{align*}
Using Plancherel's identity and noticing that $m$ is even and $u$ is $\R$-valued, we get
\[E_I^2(t)=\Lambda_2(m(\xi_1)m(\xi_2)).\]
Now, \eqref{eq:menergy} in Lemma \ref{lem51} and symmetry (about $\xi_2$ and $\xi_3$) are used to yield
\begin{align*}
\frac{d}{dt}E_I^2(t)&=\Lambda_2(m(\xi_1)m(\xi_2)v_2)-\Lambda_2(m(\xi_1)m(\xi_2)h_2)\\
&\quad-i\Lambda_3(m(\xi_1)m(\xi_2+\xi_3)(\xi_2+\xi_3))\\
&=\Lambda_3(-i[m(\xi_1)m(\xi_2+\xi_3)(\xi_2+\xi_3)]_{sym}).
\end{align*}

Putting
\begin{align*}
M_3(\xi_1,\xi_2,\xi_3)=-i[m(\xi_1)m(\xi_2+\xi_3)(\xi_2+\xi_3)]_{sym},
\end{align*}
we define the following new modified energy
\[E_I^3(t)=E_I^2(t)+\Lambda_3(\sigma_3),
\]
where the symmetric function $\sigma_3$ will be chosen momentarily
to achieve a cancelation. Applying \eqref{eq:menergy} of Lemma \ref{lem51} we get
\begin{align}\label{eq:E3}
\frac{d}{dt}E_I^3(t)&=\Lambda_3(M_3)+\Lambda_3(\sigma_3v_3)-\Lambda_3(\sigma_3h_3)\\
&-\half
3 i\Lambda_4(\sigma_3(\xi_1,\xi_2,\xi_3+\xi_4)(\xi_3+\xi_4)).\nonumber
\end{align}
Unlike the KdV case in \cite{I-method}, there is one more term to be
canceled. Thus, we choose
\begin{align*}
\sigma_3=\frac{M_3}{h_3-v_3}
\end{align*}
to force that the part containing the $\Lambda_3$ terms in \eqref{eq:E3} vanishes. So, if
\begin{align*}
M_4(\xi_1,\xi_2,\xi_3,\xi_4)=-i\half
3[\sigma_3(\xi_1,\xi_2,\xi_3+\xi_4)(\xi_3+\xi_4)]_{sym},
\end{align*}
then
\begin{align*}
\frac{d}{dt}E_I^3(t)=\Lambda_4(M_4).
\end{align*}
Similarly, if
\[E_I^4(t)=E_I^3(t)+\Lambda_4(\sigma_4)\]
with
\begin{align*}
\sigma_4=\frac{M_4}{h_4-v_4},
\end{align*}
then
\begin{align*}
\frac{d}{dt}E_I^4(t)=\Lambda_5(M_5),
\end{align*}
where
\begin{align*}
M_5(\xi_1,\ldots,\xi_5)=-2i[\sigma_4(\xi_1,\xi_2,\xi_3,\xi_4+\xi_5)(\xi_4+\xi_5)]_{sym}.
\end{align*}

In order to prove the pointwise estimates for the multipliers
$\sigma_3,\sigma_4$, we need two more lemmas.

\begin{lemma}[equalities on hyper-planes]\label{lem:decreso}
\item{\rm(i)} If $\xi_1+\xi_2+\xi_3=0$, then
\[\xi_1^3+\xi_2^3+\xi_3^3=3\xi_1\xi_2\xi_3\]
and
\[\xi_1|\xi_1|+\xi_2|\xi_2|+\xi_3|\xi_3|=2\frac{\xi_1\xi_2\xi_3}{|\xi|_{max}},\]
where $|\xi|_{max}=\max\{|\xi_j|:\,j=1,\,2,\,3\}.$

\item{\rm(ii)} If $\xi_1+\xi_2+\xi_3+\xi_4=0$, then
\[\xi_1^3+\xi_2^3+\xi_3^3+\xi_4^3=-3(\xi_1+\xi_2)(\xi_1+\xi_3)(\xi_2+\xi_3)\]
and
\[|v_4-h_4|\sim |(\xi_1+\xi_2)(\xi_1+\xi_3)(\xi_2+\xi_3)|,\]
whenever $\max\{|\xi_j|:\,j=1,\,2,\,3,\,4\}\gg 1$ and $|\tilde\alpha|\leq
1.$
\end{lemma}
\begin{proof} This follows from a straightforward computation.
\end{proof}

To introduce the next lemma, we first observe that if $m$ is of the form \eqref{eq:m} then $m^2$ enjoys
\begin{eqnarray}\label{eq:eImul}
\left\{
\begin{array}{l}
m^2(\xi)\sim m^2(\xi') \mbox{ for } |\xi|\sim|\xi'|,\\
(m^2)'(\xi)=O(\frac{m^2(\xi)}{|\xi|}),\\
(m^2)''(\xi)=O(\frac{m^2(\xi)}{|\xi|^2}).
\end{array}\right.
\end{eqnarray}
Secondly, we need two mean value formulas which follow immediately from
the fundamental theorem of calculus: $|\eta|,|\lambda|\ll |\xi|$ implies
\begin{equation}\label{eq:mvt}
|a(\xi+\eta)-a(\xi)|\lesssim |\eta|\sup_{|\xi'|\sim |\xi|}|a'(\xi')|,
\end{equation}
and the double mean value formula
\begin{equation}\label{eq:dmvt}
|a(\xi+\eta+\lambda)-a(\xi+\eta)-a(\xi+\lambda)+a(\xi)|\lesssim |\eta||\lambda|\sup_{|\xi'|\sim |\xi|}|a''(\xi')|.
\end{equation}
In applying (\ref{eq:mvt}) and (\ref{eq:dmvt}), we are required to extend the surface supported
multiplier $\sigma_3$ to the whole space as in \cite{GuoWang}. More precisely,

\begin{lemma}[extension to entire space]\label{lem53}
If $m$ is of the form \eqref{eq:m}, then for each dyadic
$\lambda\leq \mu$ there is an extension of $\sigma_3$ from the
diagonal set
\[\{(\xi_1,\xi_2,\xi_3)\in \Gamma_3(\R), |\xi_1|\sim \lambda,\quad |\xi_2|, |\xi_3|\sim \mu\}\]
to the full dyadic set
\[\{(\xi_1,\xi_2,\xi_3)\in \R^3, |\xi_1|\sim \lambda,\quad |\xi_2|, |\xi_3|\sim \mu\}\]
which satisfies
\begin{equation}\label{eq:m3}
|\partial_1^{\beta_1}\partial_2^{\beta_2}\partial_3^{\beta_3}\sigma_3(\xi_1,\xi_2,\xi_3)|\leq
C m^2(\lambda)\mu^{-2}\lambda^{-\beta_1}\mu^{-\beta_2-\beta_3}.
\end{equation}
\end{lemma}
\begin{proof} Without loss of generality, we may assume $\max(|\xi_1|,|\xi_2|,|\xi_3|)\gg 1$ (otherwise
$\sigma_3\equiv 0$). Since
\[v_3=i(\xi_1^3+\xi_2^3+\xi_3^3)=3i\xi_1\xi_2\xi_3\]
is with a size about $\lambda \mu^2$ on the hyperplane $\xi_1+\xi_2+\xi_3=0$ and since
\begin{eqnarray*}
M_3(\xi_1,\xi_2,\xi_3)&=&-i[m(\xi_1)m(\xi_2+\xi_3)(\xi_2+\xi_3)]_{sym}\\
&=&i(m^2(\xi_1)\xi_1+m^2(\xi_2)\xi_2+m^2(\xi_3)\xi_3),
\end{eqnarray*}
is valid for $\lambda \sim \mu$, we extend $\sigma_3$ by setting
\begin{equation*}
\sigma_3(\xi_1,\xi_2,\xi_3)=-\frac{i(m^2(\xi_1)\xi_1+m^2(\xi_2)\xi_2+m^2(\xi_3)\xi_3)}{3i\xi_1\xi_2\xi_3(1-\frac{2\tilde\alpha}{3|\xi|_{max}})},
\end{equation*}
and if $\lambda\ll \mu$, we extend $\sigma_3$ by setting
\begin{equation*}
\sigma_3(\xi_1,\xi_2,\xi_3)=-\frac{i(m^2(\xi_1)\xi_1+m^2(\xi_2)\xi_2-m^2(\xi_1+\xi_2)(\xi_1+\xi_2))}{3i\xi_1\xi_2\xi_3(1-\frac{2\tilde\alpha}{3|\xi|_{max}})}.
\end{equation*}
 From \eqref{eq:mvt} and \eqref{eq:eImul}, we see that \eqref{eq:m3} holds.
\end{proof}

With the foregoing treatment and some ideas in \cite{GuoWang}, we are ready to give the pointwise bounds for $\sigma_4$ which is the key to
control the growth of $E^4_I(t)$ and hence like no others (including the KdV case).

\begin{lemma}[$M_4$ estimate]\label{lem54} If $m$ is of the form \eqref{eq:m}, then
\begin{equation}\label{eq:m4}
\frac{|M_4(\xi_1,\xi_2,\xi_3,\xi_4)|}{|v_4-h_4|}\lesssim
\frac{m^2(\min(N_i,N_{jk}))}{(N+N_1)(N+N_2)(N+N_3)(N+N_4)}
\end{equation}
holds for $|\xi_i|\sim N_i,|\xi_j+\xi_k|\sim N_{jk}$ with $N_i, N_{jk}$ dyadic.
\end{lemma}

\begin{proof} From Lemma \ref{lem:decreso} it is seen that (\ref{eq:m4}) follows from a verification of

\begin{align}\label{eq:m4A}
\frac{|M_4(\xi_1,\xi_2,\xi_3,\xi_4)|}{|v_4|}\lesssim
\frac{m^2(\min(N_i,N_{jk}))}{(N+N_1)(N+N_2)(N+N_3)(N+N_4)}.
\end{align}
By symmetry, we may assume $N_1\geq N_2\geq N_3\geq N_4$.
Note that $\xi_1+\xi_2+\xi_3+\xi_4=0$ yields $N_1\sim N_2$. So, we may also
assume that $N_1\sim N_2 \ges N$ --otherwise $M_4$ vanishes due to
$m^2(\xi)=1$ when $|\xi|\leq N$.

If $\max(N_{12},N_{13},N_{14})\ll N_1$, then $\xi_2\approx-\xi_1,\ \xi_3\approx-\xi_1,\ \xi_4\approx
-\xi_1$, which contradicts that $\xi_1+\xi_2+\xi_3+\xi_4=0$. Hence
we get $\max(N_{12},N_{13},N_{14})\sim N_1$. Consequently, we rewrite the
right-hand side of \eqref{eq:m4A} as
\begin{equation*}
\frac{m^2(\min(N_i,N_{jk}))}{{N_1}^2(N+N_3)(N+N_4)}.
\end{equation*}
Using Lemma \ref{lem:decreso} we get that if $\xi_1+\xi_2+\xi_3+\xi_4=0$ then
\[v_4=i(\xi_1^3+\xi_2^3+\xi_3^3+\xi_4^3)=-3i(\xi_1+\xi_2)(\xi_1+\xi_3)(\xi_2+\xi_3)\]
is with size $N_{12}N_{13}N_{14}$. The construction of $M_4$ tells us

\begin{align}\label{eq:rm4}
M_4(\xi_1,\xi_2,\xi_3,\xi_4)\sim &\ [\sigma_3(\xi_1,\xi_2,\xi_3+\xi_4)(\xi_3+\xi_4)]_{sym}\nonumber\\
=&\ \sigma_3(\xi_1,\xi_2,\xi_3+\xi_4)(\xi_3+\xi_4)+\sigma_3(\xi_1,\xi_3,\xi_2+\xi_4)(\xi_2+\xi_4)\nonumber\\
&+\sigma_3(\xi_1,\xi_4,\xi_2+\xi_3)(\xi_2+\xi_3)+\sigma_3(\xi_2,\xi_3,\xi_1+\xi_4)(\xi_1+\xi_4)\nonumber\\
&+\sigma_3(\xi_2,\xi_4,\xi_1+\xi_3)(\xi_1+\xi_3)+\sigma_3(\xi_3,\xi_4,\xi_1+\xi_2)(\xi_1+\xi_2)\nonumber\\
=&\ [\sigma_3(\xi_1,\xi_2,\xi_3+\xi_4)-\sigma_3(-\xi_3,-\xi_4,\xi_3+\xi_4)](\xi_3+\xi_4)\nonumber\\
&+[\sigma_3(\xi_1,\xi_3,\xi_2+\xi_4)-\sigma_3(-\xi_2,-\xi_4,\xi_2+\xi_4)](\xi_2+\xi_4)\nonumber\\
&+[\sigma_3(\xi_1,\xi_4,\xi_2+\xi_3)-\sigma_3(-\xi_2,-\xi_3,\xi_2+\xi_3)](\xi_2+\xi_3)\nonumber\\
:=&\ I+II+III.
\end{align}
So, the inequality \eqref{eq:m4A} will follow from a case-by-case analysis.

(i) $|N_4|\ges \half{N}$. This case is divided into four subcases: $\rm(i)_1$; $\rm(ii)_2$; $\rm(i)_3$; $\rm(i)_4$ below.

$\rm(i)_1$ -- $N_{12}, N_{13}, N_{14}\ges N_1$. For this subcase, we just use \eqref{eq:m3} to get
\begin{align*}
\frac{|M_4(\xi_1,\xi_2,\xi_3,\xi_4)|}{|v_4|}\lesssim
\frac{m^2(N_4)}{N_1^4},
\end{align*}
which yields \eqref{eq:m4A}.

$\rm(i)_2$ -- $N_{12}\ll N_1$, $N_{13}\ges N_1$, $N_{14}\ges N_1$. Under this subcase, we are required to handle the contributions of I, II and III separately. For I, we employ \eqref{eq:m3} to derive
\begin{align*}
\frac{|I|}{|v_4|}\lesssim \frac{m^2(\min(N_4, N_{12}))}{N_1^4},
\end{align*}
which gives \eqref{eq:m4A}. For II, we rewrite
\begin{align*}
II=&\ [\sigma_3(\xi_1,\xi_3,\xi_2+\xi_4)-\sigma_3(-\xi_2,-\xi_4,\xi_2+\xi_4)](\xi_2+\xi_4)\nonumber\\
=&\ [\sigma_3(\xi_1,\xi_3,\xi_2+\xi_4)-\sigma_3(-\xi_2,\xi_3,\xi_2+\xi_4)](\xi_2+\xi_4)\nonumber\\
&+[\sigma_3(-\xi_2,\xi_3,\xi_2+\xi_4)-\sigma_3(-\xi_2,-\xi_4,\xi_2+\xi_4)](\xi_2+\xi_4)\nonumber\\
:=&\ II_1+II_2.
\end{align*}
If $N_{12}\ges N_3$, then using \eqref{eq:mvt}, \eqref{eq:m3} for
$II_1$ and using \eqref{eq:m3} for $II_2$, we find
\[\frac{|II|}{|v_4|}\ \les\ \frac{m^2(N_4)}{N_1^3 N_3}.\]
If $N_{12}\ll N_3$, using \eqref{eq:mvt}, \eqref{eq:m3} for both
$II_1$ and $II_2$, we get
\[\frac{|II|}{|v_4|}\ \les\  \frac{m^2(N_4)}{N_1^3 N_3}, \]
Adding $II_1$ and $II_2$, we reach \eqref{eq:m4A}. For III, we repeat the foregoing estimates for II, thereby obtaining \eqref{eq:m4A}.

$\rm(i)_3$ -- $N_{12}\ll N_1$, $N_{13}\ll N_1$, $N_{14}\ges N_1$. Since $N_{12}\ll N_1$, $N_{13}\ll N_1$ implies $N_1\sim N_2\sim N_3\sim N_4$, we deal with I, II, and III respectively. Regarding I, we rewrite
\begin{align*}
I=&\ [\sigma_3(\xi_1,\xi_2,\xi_3+\xi_4)-\sigma_3(-\xi_3,\xi_2,\xi_3+\xi_4)](\xi_3+\xi_4)\nonumber\\
&+[\sigma_3(-\xi_3,\xi_2,\xi_3+\xi_4)-\sigma_3(-\xi_3,-\xi_4,\xi_3+\xi_4)](\xi_3+\xi_4)\nonumber\\
:=&\ I_1+I_2.
\end{align*}
Using \eqref{eq:m3}, \eqref{eq:mvt} for both $I_1$ and $I_2$, we get
\begin{align*}
\frac{|I|}{|v_4|}\ \les\ \frac{m^2(N_{12})}{N_1^4}.
\end{align*}
thereby reaching \eqref{eq:m4A}. Regarding II, we just redo the above estimates for I to reach \eqref{eq:m4A}. Regarding III, we rewrite
\begin{align*}
III=&\ [\sigma_3(\xi_1,\xi_4,\xi_2+\xi_3)-\sigma_3(-\xi_2,-\xi_3,\xi_2+\xi_3)](\xi_2+\xi_3)\nonumber\\
=&\ \frac{1}{2}[\sigma_3(\xi_1,\xi_4,\xi_2+\xi_3)-\sigma_3(-\xi_2,-\xi_3,\xi_2+\xi_3)\nonumber\\
&-\sigma_3(-\xi_3,-\xi_2,\xi_2+\xi_3)+\sigma_3(\xi_4,\xi_1,\xi_2+\xi_3)](\xi_2+\xi_3).
\end{align*}
Using \eqref{eq:dmvt} four times, we have
\begin{align*}
\frac{|III|}{|v_4|}\ \les\ \frac{m^2(N_{1})}{N_1^4},
\end{align*}
thereby getting the desired estimate.

$\rm(i)_4$ -- $N_{12}\ll N_1$, $N_{13}\ges N_1$, $N_{14}\ll N_1$. This case is completely similar to $\rm(i)_3$. So, the details are omitted here.

(ii) $N_4\ll N/2$. In this case we have
$$
m^2(\min(N_i,N_{jk}))=1\quad\hbox{and}\quad N_{13}\sim
|\xi_1+\xi_3|=|\xi_2+\xi_4|\sim N_1.
$$
whence controlling \eqref{eq:m4A} in accordance with the following two subcases: $\rm(ii)_1$; $\rm(ii)_2$.

$\rm(ii)_1$ -- $N_1/4>N_{12}\ges N/2$. Since a combination of $N_4\ll N/2$ and $|\xi_3+\xi_4|=|\xi_1+\xi_2|\ges N/2$ implies $N_3\ges N/2$, using $|v_4|\sim N_{12}N_1^2$ we bound the six
terms in \eqref{eq:rm4} respectively, whence getting
\begin{align*}
\frac{|M_4|}{|v_4|}\ \les\ \frac{1}{N_1^2N_3N},
\end{align*}
which gives \eqref{eq:m4A}.

$\rm(ii)_2$ -- $N_{12}\ll N/2$. Owing to $N_{12}=N_{34}\ll N/2$ and $N_4\ll N/2$, we must have
$N_3\ll N/2$ and $N_{13}\sim N_{14}\sim N_1$, thereby treating I, II and III. Concerning I,  we use $N_3, N_4, N_{34}\ll N/2$ to produce $\sigma_3(-\xi_3,-\xi_4,\xi_3+\xi_4)=0$. Thus, it follows from
\eqref{eq:m3} that
\begin{align*}
\frac{|I|}{|v_4|}\ \les\
\frac{|\sigma_3(\xi_1,\xi_2,\xi_3+\xi_4)|}{N_1^2}\ \les\ \frac{1}{N_1^4},
\end{align*}
as desired. Concerning II and III, we have two items of $N_3, N_4, N_{12}$
in the denominator which will cause a problem. Thus, we cannot deal
with II and III separately, but we need to exploit the cancelation
between II and III. To do so, we rewrite
\begin{align*}
II+III=&\ [\sigma_3(\xi_1,\xi_3,\xi_2+\xi_4)-\sigma_3(-\xi_2,-\xi_4,\xi_2+\xi_4)](\xi_2+\xi_4)\nonumber\\
&+[\sigma_3(\xi_1,\xi_4,\xi_2+\xi_3)-\sigma_3(-\xi_2,-\xi_3,\xi_2+\xi_3)](\xi_2+\xi_3)\nonumber\\
=&\ [\sigma_3(\xi_1,\xi_3,\xi_2+\xi_4)-\sigma_3(-\xi_2,-\xi_4,\xi_2+\xi_4)]\xi_4\nonumber\\
&+[\sigma_3(\xi_1,\xi_4,\xi_2+\xi_3)-\sigma_3(-\xi_2,-\xi_3,\xi_2+\xi_3)]\xi_3\nonumber\\
&+[\sigma_3(\xi_1,\xi_3,\xi_2+\xi_4)-\sigma_3(-\xi_2,-\xi_4,\xi_2+\xi_4)\nonumber\\
&+\sigma_3(\xi_1,\xi_4,\xi_2+\xi_3)-\sigma_3(-\xi_2,-\xi_3,\xi_2+\xi_3)]\xi_2\nonumber\\
=&\ J_1+J_2+J_3.
\end{align*}
The consideration of $J_1$ comes first. Noticing
\begin{align*}
\frac{|J_1|}{|v_4|}\ \les\
\frac{|[\sigma_3(\xi_1,\xi_3,\xi_2+\xi_4)-\sigma_3(-\xi_2,-\xi_4,\xi_2+\xi_4)]\xi_4|}{N_{12}N_1^2},
\end{align*}
we obtain that if $N_{12}\ll N_3$ (in this case, $N_3\sim N_4$), then using
\eqref{eq:mvt} twice (otherwise using \eqref{eq:mvt} once and
\eqref{eq:m3}) one gets
\begin{align*}
\frac{|J_1|}{|v_4|}\ \les\ \frac{1}{N_1^4}.
\end{align*}
The treatment of $J_2$ is identical to that of $J_1$. Thus, it remains to handle
$J_3$. In doing so, we first assume that $N_{12}\ges N_3$. Then by the symmetry
of $\sigma_3$, we get
\begin{align*}
J_3=&\ [\sigma_3(\xi_1,\xi_3,\xi_2+\xi_4)-\sigma_3(-\xi_2-\xi_3,\xi_3,\xi_2)\nonumber\\
&+\sigma_3(\xi_1,\xi_4,\xi_2+\xi_3)-\sigma_3(-\xi_2-\xi_4,\xi_4,\xi_2)]\xi_2.
\end{align*}
From \eqref{eq:mvt} and $N_{12}\ges N_3$, we achieve
\begin{align*}
\frac{|J_3|}{|v_4|}\ \les\ \frac{1}{N_1^4}.
\end{align*}
Secondly, if $N_{12}\ll N_3$, then $N_3\sim N_4$, and hence we rewrite
\begin{align*}
J_3=&\ [\sigma_3(-\xi_2,\xi_3,\xi_2+\xi_4)-\sigma_3(-\xi_2,-\xi_4,\xi_2+\xi_4)\nonumber\\
&+\sigma_3(\xi_1,\xi_4,\xi_2+\xi_3)-\sigma_3(\xi_1,-\xi_3,\xi_2+\xi_3)]\xi_2\nonumber\\
&+[\sigma_3(\xi_1,\xi_3,\xi_2+\xi_4)-\sigma_3(-\xi_2,\xi_3,\xi_2+\xi_4)\nonumber\\
&+\sigma_3(\xi_1,-\xi_3,\xi_2+\xi_3)-\sigma_3(-\xi_2,-\xi_3,\xi_2+\xi_3)]\xi_2\nonumber\\
:=&\ J_{31}+J_{32}.
\end{align*}
On the one hand, \eqref{eq:mvt} derives
\begin{align*}
\frac{|J_{32}|}{|v_4|}\ \les\ \frac{1}{N_1^4}.
\end{align*}
On the other hand, it follows from \eqref{eq:m3} and
$m^2(\xi_3)=m^2(\xi_4)=1$ that
\begin{align*}
J_{31}&=[{\sigma}_3(-\xi_2,\xi_3,\xi_2+\xi_4)-{\sigma}_3(-\xi_2,-\xi_4,\xi_2+\xi_4)\nonumber\\
&\quad-{\sigma}_3(\xi_1,-\xi_3,\xi_2+\xi_3)+{\sigma}_3(\xi_1,\xi_4,\xi_2+\xi_3)]\xi_2\nonumber\\
&=\frac{-m^2(\xi_2)\xi_2+\xi_3+m^2(\xi_2+\xi_4)(\xi_2+\xi_4)}{\xi_3(\xi_2+\xi_4)(-\xi_2-\frac{2\tilde\alpha}{3})}\xi_2\\
&\quad-\frac{-m^2(\xi_2)\xi_2-\xi_4+m^2(\xi_2+\xi_4)(\xi_2+\xi_4)}{(-\xi_4)(\xi_2+\xi_4)(-\xi_2-\frac{2\tilde\alpha}{3})}\xi_2\nonumber\\
&\quad+\frac{m^2(\xi_1)\xi_1+\xi_4+m^2(\xi_2+\xi_3)(\xi_2+\xi_3)}{\xi_1\xi_4(\xi_2+\xi_3+\frac{2\tilde\alpha}{3})}\xi_2\\
&\quad-\frac{m^2(\xi_1)\xi_1-\xi_3+m^2(\xi_2+\xi_3)(\xi_2+\xi_3)}{\xi_1(-\xi_3)(\xi_2+\xi_3+\frac{2\tilde\alpha}{3})}\xi_2\\
&=\frac{m^2(\xi_2+\xi_4)(\xi_2+\xi_4)-m^2(\xi_2)\xi_2}{\xi_3(\xi_2+\xi_4)(-\xi_2-\frac{2\tilde\alpha}{3})}\xi_2-
\frac{m^2(\xi_2+\xi_4)(\xi_2+\xi_4)-m^2(\xi_2)\xi_2}{(-\xi_4)(\xi_2+\xi_4)(-\xi_2-\frac{2\tilde\alpha}{3})}\xi_2\\
&\quad+\frac{m^2(\xi_1)\xi_1+m^2(\xi_2+\xi_3)(\xi_2+\xi_3)}{\xi_1\xi_4(\xi_2+\xi_3+\frac{2\tilde\alpha}{3})}\xi_2-\frac{m^2(\xi_1)\xi_1+m^2(\xi_2+\xi_3)(\xi_2+\xi_3)}
{\xi_1(-\xi_3)(\xi_2+\xi_3+\frac{2\tilde\alpha}{3})}\xi_2\\
&=-\Big(\frac{\xi_3+\xi_4}{\xi_3\xi_4}\Big)\\
&\quad\times\Big[\frac{-m^2(\xi_2)\xi_2+m^2(\xi_2+\xi_4)(\xi_2+\xi_4)+m^2(\xi_1)\xi_1+m^2(\xi_2+\xi_3)(\xi_2+\xi_3)}{(\xi_2+\xi_4)(\xi_2+\frac{2\tilde\alpha}{3})}\Big]\xi_2\\
&\quad+\Big(\frac{\xi_3+\xi_4}{\xi_3\xi_4}\Big)[m^2(\xi_1)\xi_1+m^2(\xi_2+\xi_3)(\xi_2+\xi_3)]\\
&\quad\times\Big[\frac{1}{\xi_1(\xi_2+\xi_3+\frac{2\tilde\alpha}{3})}+\frac{1}{(\xi_2+\xi_4)(\xi_2+\frac{2\tilde\alpha}{3})}\Big]\xi_2.
\end{align*}
Therefore, we use \eqref{eq:dmvt} for the first term, and
\eqref{eq:mvt} for the second term, to conclude
\begin{align*}
\frac{|J_{31}|}{|v_4|}\ \les\ \frac{1}{N_1^4},
\end{align*}
which completes the estimate of $J_3$.
\end{proof}

Below is the estimate for $M_5(\xi_1,...,x_5)$.

\begin{lemma}[$M_5$ estimate]\label{prop:lowchM5} If $m$ is of the form \eqref{eq:m}, then
\begin{align*}
|M_5(\xi_1,\ldots,\xi_5)|\ \les\
\left[\frac{m^2(N_{*45})N_{45}}{(N+N_1)(N+N_2)(N+N_3)(N+N_{45})}\right]_{sym},
\end{align*}
where
\[N_{*45}=\min(N_1,N_2,N_3,N_{45},N_{12},N_{13},N_{23}).\]
\end{lemma}
\begin{proof} This is immediate from the estimates of $\sigma_4$ in Lemma \ref{lem54}.
\end{proof}

\section{Proof of Theorem \ref{thm4}}

Through demonstrating Theorem \ref{thm4}, we, in this section, extend the local solutions in Theorems \ref{thm2}-\ref{thm3} to the global solutions.
The argument depends on a variant of the local well-posedness as follows.

\begin{proposition}[variant of local well-posedness]\label{prop:lowchvlwp}
Let $-3/4\leq s\leq 0$. Suppose $\phi$ satisfies
$\norm{I\phi}_{L^2(\R)}\leq 2\epsilon_0\ll 1$. Then
\eqref{sBen} has a unique solution on $[-1,1]$ with
\begin{eqnarray*}
\norm{Iu}_{\bar{F}^s(1)}\leq C\epsilon_0,
\end{eqnarray*}
where $C$ is a positive constant independent of $N$ and $0<\lambda\leq 1$.
\end{proposition}
\begin{proof} This can be verified via a slight modification of that for Theorem \ref{thm3}.
\end{proof}

From Proposition \ref{prop:lowchvlwp}, we see that it is enough to control
the growth of $E_I^2(t)$. In doing so, it is better to settle directly the
growth of $E_I^4(t)$ via the following difference inequality at the intermediate point $s=-3/4$.

\begin{proposition}[difference between $E^4_I$ and $E^2_I$]\label{prop:EI2EI4}
Let $s=-3/4$ and $I$ be defined with the multiplier $m$ of the form \eqref{eq:m}. Then
\begin{equation*}
|E_I^4(t)-E_I^2(t)|\ \les\ \norm{Iu(t)}_{L^2}^3+\norm{Iu(t)}_{L^2}^4.
\end{equation*}
\end{proposition}

\begin{proof} Because of
$$
E_I^4(t)=E_I^2(t)+\Lambda_3(\sigma_3)+\Lambda_4(\sigma_4),
$$
it suffices to show the following two inequalities:
\begin{eqnarray}
|\Lambda_3(\sigma_3;u_1,u_2,u_3)|&\les&
\prod_{i=1}^3\norm{Iu_i}_{L^2};\label{eq:endiff3linear}\\
|\Lambda_4(\sigma_4;u_1,u_2,u_3,u_4)|&\les&
\prod_{i=1}^4\norm{Iu_i}_{L^2}.\label{eq:endiff4linear}
\end{eqnarray}
In the sequel, we may assume that $\wh{u_i}$ are non-negative.

To prove \eqref{eq:endiff3linear}, it suffices to check
\begin{align}\label{eq:lowchE3linearprf1}
\left|{\Lambda_3\brk{\frac{m^2(\xi_1)\xi_1+m^2(\xi_2)\xi_2+m^2(\xi_3)\xi_3}{\xi_1\xi_2\xi_3m(\xi_1)m(\xi_2)m(\xi_3)};u_1,u_2,u_3}}\right|\ \les\ \prod_{i=1}^3\norm{u_i}_2.
\end{align}
By the Littlewood-Paley decomposition, we find that the left-hand side of
\eqref{eq:lowchE3linearprf1} is bounded by
\begin{align*}\label{eq:lowchE3linearprf2}
\sum_{\Lambda_3}:=\sum_{k_i\geq
0}\left|{\Lambda_3\brk{\frac{m^2(\xi_1)\xi_1+m^2(\xi_2)\xi_2+m^2(\xi_3)\xi_3}{\xi_1\xi_2\xi_3m(\xi_1)m(\xi_2)m(\xi_3)};P_{k_1}u_1,P_{k_2}u_2,P_{k_3}u_3}}\right|.
\end{align*}
Let $N_i=2^{k_i}$. Using symmetry we may also assume $N_1\geq N_2 \geq
N_3$ and hence $N_1\sim N_2\ges N$. Consequently, we need to handle two cases.

(i) $N_3\ll N$. This case ensures $m(N_3)=1$, but also
\begin{eqnarray*}
\sum_{\Lambda_3}&\les&\sum_{k_i\geq
0}\left|{\Lambda_3\brk{\frac{N^sN^s}{N_1^{1+s}N_2^{1+s}};P_{k_1}u_1,P_{k_2}u_2,P_{k_3}u_3}}\right|\\
&\les&\sum_{k_i\geq
0}\left|{\Lambda_3\brk{N_1^{-1/4}N_2^{-1/4};P_{k_1}u_1,P_{k_2}u_2,P_{k_3}u_3}}\right|.
\end{eqnarray*}
So, (\ref{eq:lowchE3linearprf1}) will be proved upon verifying
\begin{equation*}
\sum_{k_i\geq 0}\int_{\xi_1+\xi_2+\xi_3=0,|\xi_i|\sim
N_i}N_1^{-1/2}\prod_{i=1}^3\eta_{k_i}(\xi_i)\wh{u_i}(\xi_i)\ \les\ \prod_{i=1}^3\norm{u_i}_{L^2}.
\end{equation*}
To see this, let us define $v_i(x)$ via its Fourier transform:
\[\wh{v_i}(\xi)=N_i^{-1/6}\wh{u_i}(\xi)\chi_{\{|\xi|\sim N_i\}}(\xi).\]
By the Sobolev embedding inequality we have $\norm{v_i}_{L^3}\les
\norm{u_i}_{L^2}$, thus getting by H\"older's inequality,
\begin{eqnarray*}
\sum_{k_i\geq 0}\int_{\xi_1+\xi_2+\xi_3=0,|\xi_i|\sim
N_i}N_1^{-1/2}\prod_{i=1}^3\eta_{k_i}(\xi_i)\wh{u_i}(\xi_i)&\les&\sum_{k_i\geq
0}N_1^{-1/6}N_3^{1/6}\prod_{i=1}^3\norm{v_i}_{L^3}\\
&\les&\prod_{i=1}^3\norm{u_i}_{L^2},
\end{eqnarray*}
as desired.

(ii) $N_3\ges N$. Under this assumption, it is not hard to obtain
\begin{align*}
\sum_{\Lambda_3}\ \les\ \sum_{k_i\geq
0}\left|{\Lambda_3\brk{\frac{N_3^{-3/4}N^{-3/4}}{N_1^{1/2}};P_{k_1}u_1,P_{k_2}u_2,P_{k_3}u_3}}\right|\ \les\ \prod_{i=1}^3\norm{u_i}_{L^2},
\end{align*}
whence getting \eqref{eq:endiff3linear}.

Next, in order to show \eqref{eq:endiff4linear}, it is enough to prove
\begin{eqnarray}\label{eq:lowchE4linearprf1}
\left|{\Lambda_4\brk{\frac{\sigma_4}{m(\xi_1)m(\xi_2)m(\xi_3)m(\xi_4)};u_1,u_2,u_3,u_4}}\right|\ \les\
\prod_{i=1}^4\norm{u_i}_2.
\end{eqnarray}
Again, by the Littlewood-Paley decomposition we find that the left-hand side of
\eqref{eq:lowchE4linearprf1} is dominated by
\begin{eqnarray*}\label{eq:lowchE4linearprf2}
\sum_{\Lambda_4}:=\sum_{k_i\geq
0}\left|{\Lambda_4\brk{\frac{\sigma_4}{m(\xi_1)m(\xi_2)m(\xi_3)m(\xi_4)};P_{k_1}u_1,P_{k_2}u_2,P_{k_3}u_3,P_{k_4}u_4}}\right|.
\end{eqnarray*}
Let $N_i=2^{k_i}$. Using symmetry we may assume $N_1\geq N_2 \geq
N_3\geq N_4$ and hence $N_1\sim N_2\ges N$. Thanks to
\[\left|\frac{\sigma_4}{m(\xi_1)m(\xi_2)m(\xi_3)m(\xi_4)}\right|\ \les\ \frac{1}{\prod_{i=1}^4(N+N_i)m(N_i)}\ \les\ \frac{N^{-3}}{\prod_{i=1}^4N_i^{1/4}},\]
using H\"older's inequality we get
\begin{eqnarray*}
\sum_{\Lambda_4}&\les& \sum_{k_i\geq
0}\frac{N^{-3}}{\prod_{i=1}^4N_i^{1/4}}\norm{P_{k_1}u_1}_{L^2}\norm{P_{k_2}u_2}_{L^2}\norm{P_{k_3}u_3}_{L^\infty}\norm{P_{k_4}u_4}_{L^\infty}\\
&\les&\prod_{i=1}^4\norm{u_i}_2,
\end{eqnarray*}
as desired.
\end{proof}

According to Proposition \ref{prop:EI2EI4}, $E_I^2(t)$ is very close to $E_I^4(t)$, so our task is in turn to control
$E_I^4(t)$. In order to handle the
increasing of $E_I^4(t)$, we induce the forthcoming product estimate to control the derivative
\[\frac{d}{dt}E_I^4(t)=\Lambda_5(M_5),\]
where
\[M_5(\xi_1,\ldots,\xi_5)=-2i[\sigma_4(\xi_1,\xi_2,\xi_3,\xi_4+\xi_5)(\xi_4+\xi_5)]_{sym}.\]

\begin{proposition}[product estimate]\label{pro}
Let $I\subset \R$ with $|I|\les 1$. If $0\leq k_1\leq \ldots
\leq k_5$ and $k_4\geq 10$, then
\begin{equation}\label{eq:5linear1}
\left|\int_I\int \prod_{i=1}^5 P_{k_i}(w_i)(x,t)dxdt\right|\ \les\
2^{\frac{5}{12}(k_1+k_2+k_3)}2^{-k_4-k_5}\prod_{j=1}^5\norm{\wh{P_{k_j}(w_j)}}_{X_{k_j}},
\end{equation}
where $X_{k_j}$ is replaced by $\bar{X}_{k_j}$ on
the right-hand side whenever $k_j=0$.
\end{proposition}
\begin{proof}
From H\"older's inequality it follows that the left-hand side of \eqref{eq:5linear1}
is dominated by
\[\prod_{i=1}^3\norm{P_{k_i}(w_i)}_{L_x^3L_{t\in I}^\infty}\cdot \norm{P_{k_4}(w_4)}_{L_x^\infty L_{t}^2}\cdot \norm{P_{k_5}(w_5)}_{L_x^\infty L_{t}^2}.\]
Then we use Lemma \ref{lem2} to dominate $\norm{P_{k_4}(w_4)}_{L_x^\infty L_{t}^2}$ and
$\norm{P_{k_5}(w_5)}_{L_x^\infty L_{t}^2}$. However, in controlling $\norm{P_{k_3}(w_3)}_{L_x^3L_{t\in I}^\infty}$ we consider two cases: First, if $k_3\geq 10$ then
we use interpolation between $\norm{P_{k_3}(w_3)}_{L_x^2L_{t\in
I}^\infty}$ and $\norm{P_{k_3}(w_3)}_{L_x^4L_{t\in I}^\infty}$,
Lemma \ref{lem2}. Second, if $k_3\leq 10$, then we use interpolation between
$\norm{P_{k_3}(w_3)}_{L_x^2L_{t\in I}^\infty}$ and
$\norm{P_{k_3}(w_3)}_{L_x^\infty L_{t\in I}^\infty}$. Similarly, we can handle
the remaining items: $\norm{P_{k_1}(w_1)}_{L_x^\infty L_{t}^2}$ and
$\norm{P_{k_2}(w_2)}_{L_x^\infty L_{t}^2}$, thereby reaching (\ref{eq:5linear1}).
\end{proof}

The following is the required integral inequality which has a root in \cite[Lemma 5.2]{I-method}.

\begin{proposition}[integral estimate for $\Lambda_5(M_5;\cdots)$]\label{prop:inEI4}
Let $\delta\les 1$. If $m$ is of the form \eqref{eq:m} with
$s=-3/4$, then
\begin{eqnarray}
\left|\int_0^\delta \Lambda_5(M_5;u_1,\ldots,u_5)dt\right|\ \les\
N^{-\frac{15}{4}} \prod_{j=1}^5\norm{Iu_j}_{\bar F^{0}(\delta)}.
\end{eqnarray}
\end{proposition}
\begin{proof} Without loss of generality, we may assume each $\hat{u_j}$ is nonnegative, and then restrict it to a frequency band $|\xi_j|\sim N_j=2^{k_j}$ via a Littlewood-Paley decomposition, and finally sum in $N_j$. According to Lemma \ref{lem54} and the definition of the operator $I$,  it is sufficient for us to
demonstrate
\begin{equation}\label{eqEq}
\Big|\int_0^\delta\Big(\frac{N_{45}m^2(N_{\ast 45})\prod_{j=1}^5 [m(N_j)]^{-1}}{(N+N_{45})\prod_{j=1}^3(N+N_j)}; u_1,...,u_5\Big)\,dt\Big|\lesssim N^{-\frac{15}{4}}\prod_{j=1}^5\|u_j\|_{\bar{F}^0(\delta)}.
\end{equation}
Upon canceling $N_{45}\le(N+N_{45})$ in (\ref{eqEq}), we may only deal with the worst situation when $m^2(N_{\ast 45})=1$ in the sequel, and are about to verify the following inequality

\begin{eqnarray*}
&&\sum_{k_1,\ldots, k_5\geq 0}\left|\int_0^\delta
\Lambda_5\left(\frac{\prod_{i=1}^3{[(N+N_i)m(N_i)]^{-1}}}{m(N_4)m(N_5)};P_{k_1}u_1,\ldots,P_{k_5}u_5\right)dt\right|\\
&&\ \les\ N^{-\frac{15}{4}}\prod_{i=1}^5
\norm{u_i}_{\bar{F}^0(\delta)},\quad\hbox{where}\quad N_i=2^{k_i}.
\end{eqnarray*}
Applying symmetry, we may assume $N_1\geq N_2\geq
N_3$ and $N_4\geq N_5$ and two of the estimates $N_i\ges N$. We fix the
extension $\wt{u}_i$, still denoted by $u_i$, such that $\norm{\wt{u}_i}_{\bar{F}^0}\les
2\norm{{u}_i}_{\bar{F}^0(\delta)}$.

Using the form \eqref{eq:m} with $s=-3/4$ we get
$$
\frac{1}{(N+N_i)m(N_i)}\ \les\ N^{-3/4}\jb{N_i}^{-1/4}\quad\hbox{and}\quad
\frac{1}{m(N_4)m(N_5)}\ \les\ N^{-3/2}N_4^{3/4}N_5^{3/4}.
$$
Therefore, we need to control
\begin{align*}\label{eq:deriE4}
\sum_N:=N^{-\frac{15}{4}}\sum_{k_i}\int_0^\delta
\Lambda_5\left({\jb{N_1}^{-1/4}\jb{N_2}^{-1/4}\jb{N_3}^{-1/4}N_4^{3/4}N_5^{3/4}};u_1,\ldots,u_5\right)dt.
\end{align*}
If $N_2\sim N_1 \ges N$, $N_4\les N_2$, then we just consider the worst case
$$N_1\geq N_2\geq N_4\geq N_5\geq N_3.$$
From \eqref{eq:5linear1} in Proposition \ref{pro} we see
\begin{align*}
\sum_N\ \les\ & N^{-\frac{15}{4}}\sum_{N_i}
\jb{N_1}^{-5/4}\jb{N_2}^{-5/4}\jb{N_3}^{1/6}N_4^{7/6}N_5^{7/6}\prod_{i=1}^5
\norm{\wh{P_{k_i}u}}_{X_{k_i}}\nonumber\\
\les\ & N^{-\frac{15}{4}} \prod_{j=1}^5\norm{Iu_j}_{\bar
F^{0}(\delta)}.
\end{align*}
The rest cases: $N_4\sim N_5 \ges N$, $N_1\les N_5$ or $N_1\sim N_4
\ges N$ follow in a similar ways, and so their details are omitted here.
\end{proof}

With the previous propositions, we can now extend the local solutions to the global ones, thereby completing the proof of Theorem \ref{thm4}. To see this, let us
fix $u_0\in H^{s}$ and time $T>0$. Then, our goal is to
construct the solution of \eqref{Ben1} on $t\in [0,T]$. If $u$ is a
solution to \eqref{Ben1} with initial data $u_0$, then for any
$\lambda>0$, $u_\lambda(x,t)=\lambda^{2}u(\lambda x, \lambda^3 t)$
is a solution to \eqref{sBen} with initial data
$u_{0,\lambda}(x)=\lambda^{2}u_0(\lambda x)$. By a simple calculation we
know that for $s\ge -3/4$,
\[
\norm{Iu_{0,\lambda}}_{L^2}\ \les\ \lambda^{\frac{3}{2}+s}N^{-s}\norm{u_0}_{H^s}.
\]
For a fixed $N$ (to be determined later) and $\phi(x)=\lambda^2u_0(\lambda x)$, we take $\lambda\sim
N^{\frac{2s}{3+2s}}$ such that
\begin{align*}
\lambda^{\frac32+s}N^{-s}\norm{\phi}_{H^{s}}=\epsilon_0<1.
\end{align*}
So, the goal is in turn to find the solution to
\eqref{sBen} on $[0,\lambda^{-3} T]$. According to Proposition
\ref{prop:lowchvlwp} we have a local solution $u_\lambda$ on $t\in [0,1]$,
so what we need is to control the modified energy $E_I^2(t):=\norm{Iu_\lambda}^2_{L^2}$.

Let us start with controling $E_I^2(t)$ for $t\in [0,1]$. This will be done via
proving $E_I^2(t)<4\epsilon_0^2$. Using the so-called bootstrap argument we may
assume $E_I^2(t)<5\epsilon_0^2$. From Proposition
\ref{prop:EI2EI4} we get
\[E_I^4(0)=E_I^2(0)+O(\epsilon_0^3)\quad\hbox{and}\quad E_I^4(t)=E_I^2(t)+O(\epsilon_0^3),\]
thereby finding from Proposition \ref{prop:inEI4} that for $t\in [0,1]$
\[E_I^4(t)\leq E_I^4(0)+C\epsilon_0^5N^{-15/4}.\]
Therefore
\begin{align*}
\norm{Iu_\lambda(1)}_{L^2}^2=&\ E_I^4(1)+O(\epsilon_0^3)\\
\leq&\ E_I^4(0)+C\epsilon_0^5N^{-35/4}+O(\epsilon_0^3)\\
=&\ \epsilon_0^2+C\epsilon_0^5N^{-15/4}+O(\epsilon_0^3)\\
<&\ 4\epsilon_0^2.
\end{align*}
Consequently, $u_\lambda$ is extendable to the interval $[0,2]$. Continuing this process
$M$-steps, we get that for $t\in [0,M+1]$ there is a constant $C>0$ with
\begin{align*}
E_I^4(t)\leq E_I^4(0)+CM\epsilon_0^5N^{-15/4}.
\end{align*}
Now, as long as $MN^{-15/4}\les 1$, we have
\[E_I^2(M)=E_I^4(t)+O(\epsilon_0^3)=\epsilon_0^2+O(\epsilon_0^3)+CM\epsilon_0^5N^{-15/4}<4\epsilon_0^2,\]
Thus, the solution can be extended to the interval $[0,N^{15/4}]$. Taking $N(T)$ sufficiently large such that
\[
N^{15/4}>\lambda^{-3} T\sim N^{-\frac{6s}{3+2s}}T\quad\hbox{where}\quad s\ge -3/4.
\]
Therefore, $u_\lambda$ can be extended to $[0,\lambda^{-3} T]$. Returning to the original
equation \eqref{Ben1}, we use the scaling argument to verify that its solution $u$ can be
extended to $[0,T]$, as desired.

Last but not least, let us have a look at the selection of $N$ for the given time $T>0$. Using the scaling we get
\begin{align*}
\sup_{t\in [0,T]}\norm{u(t)}_{H^{s}}\sim&
\lambda^{{-\frac32-s}}\sup_{t\in
[0,\lambda^{-3}T]}\norm{u_\lambda(t)}_{H^{s}}\leq
\lambda^{-\frac32-s}\sup_{t\in
[0,\lambda^{-3}T]}\norm{Iu_\lambda(t)}_{L^2}
\end{align*}
and
\begin{align*}
\norm{I\phi_\lambda}_{L^2}\les
N^{-s}\norm{\phi_\lambda}_{H^{s}}\sim
N^{-s}\lambda^{\frac32+s}\norm{\phi}_{H^{s}}.
\end{align*}
From the above we know
\[\sup_{t\in [0,\lambda^{-3} T]}\norm{Iu_\lambda(t)}_{L^2}\ \les\ \norm{I\phi_\lambda}_{L^2},\]
thereby finding
\[\sup_{t\in [0,T]}\norm{u(t)}_{H^{s}}\ \les\ N^{{-s}}\norm{\phi}_{H^{s}}.\]
Taking $\lambda$ such that $\norm{I\phi_\lambda}_{L^2}\sim \epsilon_0
\ll 1$, we obtain
\[\lambda=\lambda(N,\epsilon_0,\norm{\phi}_{H^{s}})\sim
\brk{\frac{\norm{\phi}_{H^{s}}}{\epsilon_0}}^{-\frac{2}{3+2s}}N^{\frac{2s}{3+2s}}\]
and consequently choose such an $N$ that
$$
N^{\frac{15}{4}} > \lambda^{-3}
T\sim C_{\norm{\phi}_{H^s},\epsilon_0}N^{-\frac{6s}{3+2s}}T.
$$
Of course, this last requirement will be met as long as
$$
N\sim c_{\norm{\phi}_{H^s},\epsilon_0} T^{\frac{4(3+2s)}{45+54s}}
$$
holds for some positive constant $c_{\norm{\phi}_{H^s},\epsilon_0}$. In particular, if $N\sim T^{4/3}$, then the found global solution $u(x,t)$ enjoy the following time-dependent estimate
$$
\|u(\cdot,t)\|_{H^{-3/4}}\lesssim (1+|t|)\|u_0\|_{H^{-3/4}}\quad\hbox{for}\quad t\in [0,T].
$$

\medskip

\noindent{\bf Acknowledement}. This paper was completed during the third-named author's visit to Polytechnic Institute of NYU in Fall 2009. This author would like to thank the math departmental faculty and staff of the institute for their kind hospitality, but also Professor Louis Nirenberg from Courant Institute of Mathematical Sciences at NYU for his helpful encouragement.

\bibliographystyle{amsplain}

\end{document}